\documentclass[reqno]{amsart}
\usepackage[T1]{fontenc}
\usepackage{thmtools}
\usepackage{xcolor}
\usepackage[breaklinks,colorlinks,linkcolor=cyan!30!blue,citecolor=cyan!30!blue]{hyperref}
\usepackage{amssymb,amsmath}
\usepackage{amsthm}
\usepackage{graphicx}
\usepackage{bibentry}

\usepackage{comment}
\usepackage[shortlabels]{enumitem}
\definecolor{nottbeaublue}{rgb}{0.74, 0.83, 0.9}
\definecolor{nottblue}{RGB}{32, 122, 178}
\definecolor{nottblue2}{RGB}{24, 71, 128}
\newcommand{\rednote}[1]{\textcolor[hsb]{.0,1.0,.75}{#1}}

\makeatletter
\def\widebreve{\mathpalette\wide@breve}
\def\wide@breve#1#2{\sbox\z@{$#1#2$}%
     \mathop{\vbox{\m@th\ialign{##\crcr
\kern0.08em\brevefill#1{0.8\wd\z@}\crcr\noalign{\nointerlineskip}%
                    $\hss#1#2\hss$\crcr}}}\limits}
\def\brevefill#1#2{$\m@th\sbox\tw@{$#1($}%
  \hss\resizebox{#2}{\wd\tw@}{\rotatebox[origin=c]{90}{\upshape(}}\hss$}
\makeatletter

\newcommand{\norm}[2][y]{\if#1y\left\fi\lVert#2\if#1y\right\fi\rVert} 
\newcommand{\abs}[2][y]{\if#1y\left\fi\lvert#2\if#1y\right\fi\rvert} 
\newcommand{\LL}{\mathrm{L}}

\newcommand{\CC}{\mathrm{C}}
\newcommand{\PP}{\mathrm{P}}
\newcommand{\bd}[1]{\boldsymbol{#1}}
\newcommand{\tri}{\mathrm{T}}
\newcommand{\Gam}{\bd{\gamma}}

\newcommand{\Kap}{\bd{\kappa}}
\newcommand{\W}{\mathrm{W}}
\newcommand{\J}{\mathrm{J}}
\newcommand{\dd}[1][y]{\if#1y\,\fi{\mathrm d}} 
\DeclareMathOperator{\proj}{proj}

\DeclareMathOperator{\vspan}{span}

\numberwithin{equation}{section}

\declaretheorem[name=Theorem,numberwithin=section]{thm}
\declaretheorem[name=Lemma,sibling=thm]{lemma}
\declaretheorem[name=Proposition,sibling=thm]{proposition}
\declaretheorem[name=Corollary,sibling=thm]{corollary}
\declaretheorem[name=Remark,sibling=thm,style=remark]{remark}
\declaretheorem[name=Definition,sibling=thm,style=definition]{definition}

\title[Markov inequality on the simplex]{Jacobi-weighted $\boldsymbol{\LL^{2}}$ Markov inequality on the simplex}
\author{Leonardo E. Figueroa}
\email{lefiguer@udec.cl}
\address{CI\textsuperscript{2}MA and Departamento de Ingenier\'ia Matem\'atica, Universidad de Concepci\'on, Casilla 160-C, Concepci\'on, Chile}
\author{Brayan Sandoval-Le\'on}
\email{bsandoval2018@udec.cl}
\address{CI\textsuperscript{2}MA and Departamento de Ingenier\'ia Matem\'atica, Universidad de Concepci\'on, Casilla 160-C, Concepci\'on, Chile}

\begin{document}
\begin{abstract}
We establish sharp Markov- and Schur-type inequalities in $\LL^{2}$ norms with Jacobi weights in the $d$-simplex.
To this end, we exploit a Bernstein-type inequality established by Ge and Xu (Constr.\ Approx.\ 62:305--328, 2025) and certain structural properties of orthogonal polynomials.
From among the latter, we make profuse use of the covariance between Jacobi weight parameters and affine maps that permute the vertices of the simplex.

\vspace{0.3cm}
\end{abstract}
\maketitle
\section{introduction}
Let $\tri^{d}$ be the simplex defined by $\tri^{d} := \{x \in \mathbb{R}^{d} \mid x_{1},\dots,x_{d} > 0,~ 1-\abs{x}> 0 \}$, where $\abs{x} := x_{1} + \cdots + x_{d}$.
Given $\Gam \in (-1,\infty)^{d+1}$, let the Jacobi-type weight function $\W_{\Gam} \colon \tri^d \to \mathbb{R}$ be defined by
\begin{equation}\label{weight}
\W_{\Gam}(x) := x^{\gamma_{1}}_{1}\cdots x^{\gamma_{d}}_{d}(1-\abs{x})^{\gamma_{d+1}}
\end{equation}
and let $\LL^2_{\Gam}$ be the weighted Lebesgue space $\LL^2(\tri^{d},\W_{\Gam}) := \{ \W_{\Gam}^{-1/2} f \mid f \in \LL^{2}(\tri^{d}) \}$, whose natural squared norm is
\begin{equation*}
\norm{u}_{\Gam}^{2} := \int_{\tri^{d}} \abs{u(x)}^{2} \W_{\Gam}(x) \dd x.
\end{equation*}
Given an integer $n \geq 0$, let $\Pi_{n}^{d}$ be the space of $d$-variate polynomials of total degree no higher than $n$ (we adopt the convention $\Pi_{n}^{d} = \{0\}$ for $n < 0$).
The main purpose of this work is proving the following Markov-type inequality on the simplex:
For all parameters $\Gam \in (-1,\infty)^{d+1}$, there exists a constant $C > 0$ such that,
\begin{equation}\label{crude-Markov}
(\forall\,n\in\mathbb{N}_0)\ (\forall\,p\in\Pi^{d}_{n}) \quad \norm{\nabla p}_{\Gam} \leq C n^2 \norm{p}_{\Gam}.
\end{equation}
We will also show that \eqref{crude-Markov} is sharp, in the sense that, in general, the power on $n$ cannot be lowered.
As it has been the case elsewhere in the literature, we come by our Markov-type inequality with the assistance of Bernstein- and Schur-type inequalities.

Among the latter, we prove three families, of which the first two come very naturally:
For all $\Gam \in (-1,\infty)^{d+1}$ and every $i \in [d+1] := \{1, \dotsc, d+1\}$, there exist a constant $C > 0$ and an integer $n_0 \in \mathbb{N}$ such that
\begin{equation}\label{crude-natural-Schur}
(\forall\,n\geq n_0)\ (\forall\,p\in\Pi^{d}_{n}) \quad
\norm{p}_{\Gam} \leq C n \norm{p}_{\Gam + e_i}
\end{equation}
and, again for all $\Gam \in (-1,\infty)^{d+1}$ and every $i, j \in [d+1]$ with $i \neq j$, there exist a constant $C > 0$ and an integer $n_0 \in \mathbb{N}$ such that
\begin{equation}\label{crude-middle-Schur}
(\forall\,n\geq n_0)\ (\forall\,p\in\Pi^{d}_{n}) \quad
\norm{p}_{\Gam} \leq C n^2 \norm{p}_{\Gam + e_i + e_j}.
\end{equation}
The third family is:
For all $\Gam \in (-1,\infty)^{d+1}$ and every $i, j \in [d+1]$ with $i \neq j$, there exist a constant $C > 0$ and an integer $n_0 \in \mathbb{N}$ such that
\begin{equation}\label{crude-depleted-Schur}
(\forall\,n\geq n_0)\ (\forall\,p\in\Pi^{d}_{n}) \quad
\norm{p}_{\Gam} \leq C \, n \norm{\frac{p}{\sqrt{x_i+x_j}}}_{\Gam + e_i + e_j},
\end{equation}
where $x_{d+1} = 1 - \abs{x}$.
We establish each of the Schur-type inequalities \eqref{crude-natural-Schur}, \eqref{crude-middle-Schur} and \eqref{crude-depleted-Schur} by first proving a variant, with a lower power on $n$ in each case, in which only orthogonal polynomials with respect to the inner product that induces the right-hand side norm are considered.

This work contributes to the extensive literature on Markov-type inequalities akin to \eqref{crude-Markov}.
For the univariate setting on the interval $(-1,1)$ involving general families of weights, we highlight the contributions of \cite{Mastroianni-Totik:2000} and \cite{Goetgheluck:1979}, from which our result for $d=1$ can be deduced via an appropriate change of variables.
In the context of higher dimensions and general domains, the works \cite{Kroo:2009} and \cite{Goetgheluck:1987} establish analogues of \eqref{crude-Markov} for the unweighted case on convex sets and locally Lipschitz compact domains, respectively; since the simplex $\tri^{d}$ is convex, our result for $\Gam = 0$ follows immediately.
On the other hand, closely related results have been established on the unit ball in \cite{Figueroa:2017b}, \cite{Benavides:2021} and \cite{LiWangWang:2023}.
Specifically, the first two works establish bounds analogous to \eqref{crude-Markov} for Gegenbauer weights (cf.~\cite[Lem.~3.4]{Figueroa:2017b}) and generalized Gegenbauer weights (cf.~\cite[Prop.~4.8]{Benavides:2021}), respectively, whereas the third extends these findings to the broader class of doubling weights.
Finally, two other relevant contributions are the works \cite{Ditzian:1992} and \cite{Dai:2006}; the former establishes unweighted Markov-type inequalities for certain directional derivatives, while the latter derives weighted Bernstein--Markov inequalities on the sphere.

Regarding the Schur-type inequalities \eqref{crude-natural-Schur}, \eqref{crude-middle-Schur} and \eqref{crude-depleted-Schur}, analogous results have been established in \cite{Figueroa:2017b} and \cite{Benavides:2021}, where Schur-type inequalities on the unit ball for Gegenbauer and generalized Gegenbauer weights, respectively, were obtained.
It is worth noting that, in both of these papers, these inequalities exhibit a shift in only one parameter of the dominant norm and are directly compatible with Bernstein-type inequalities in order to produce sharp Markov-type inequalities.
This contrasts with the present case of the simplex because, firstly, the one-parameter shift Schur-type inequality \eqref{crude-natural-Schur} cannot be chained with a Bernstein-type inequality because the latter naturally involve two-parameter shifts; secondly, the two-parameter shift Schur-type inequality \eqref{crude-middle-Schur} can indeed be chained with a classic Bernstein-type inequality, but the high power on $n$ in \eqref{crude-middle-Schur} makes the resulting Markov-type inequality not sharp; thirdly, the two-parameter Schur-type inequality \eqref{crude-depleted-Schur} is harder to prove yet does get successfully chained to a Bernstein-type inequality essentially found in \cite{GeXu:2024} in order to produce a sharp Markov-type inequality.

We expect our Schur- and Markov-type inequalities to be relevant to the Polynomial Approximation Theory of the simplex just like their analogues on the ball are in \cite{Figueroa:2017b} and \cite{Benavides:2021}.

The outline of this article is as follows.
In \autoref{sec:Jacpoly}, we present several properties of Jacobi polynomials that will be useful in subsequent sections.
In \autoref{sec:Orthogonalpolyspace}, we introduce the orthogonal polynomial spaces on the simplex and establish some of their fundamental properties.
In \autoref{sec:SGAonOPS}, we examine the Symmetric Group action on some of the orthogonal polynomial spaces introduced in \autoref{sec:Orthogonalpolyspace} and investigate its properties.
In \autoref{sec:Schur-type-inequalities}, we prove the Schur-type inequalities \eqref{crude-natural-Schur}, \eqref{crude-middle-Schur}, and \eqref{crude-depleted-Schur}: the first is obtained via integration by parts, the second follows by combining the first inequality with Parseval's identity, and the third is established through a distinct approach relying essentially on an induction argument on the dimension of the simplex.
In \autoref{sec:markovinequality}, we prove the Markov-type inequality \eqref{crude-Markov}.
The proof strategy for the main result therein, \autoref{thm:Markovineq}, comprises two stages: first, we establish the inequality for the directional derivative along $e_{1}-e_{2}$ by combining a Schur-type inequality (cf.~\autoref{thm:whole-depleted-Schur}) with a Bernstein-type directional derivative bound (cf.~\autoref{prop:Bernstein12}); second, we extend the result to the canonical directional derivatives with the help of the results of \autoref{sec:SGAonOPS}.
\section{Jacobi Polynomials}\label{sec:Jacpoly}
Given $\alpha,\beta>-1$ let the weight function $\J(\alpha,\beta):(-1,1)\longrightarrow\mathbb{R}$ be defined by
\begin{equation}\label{eq:Jacobiweight}
\J(\alpha,\beta)(t) = (1-t)^{\alpha}(1+t)^{\beta}.
\end{equation}
The Jacobi polynomials $\PP_{n}^{(\alpha,\beta)}$ are the orthogonal polynomials on $(-1,1)$ with respect to the inner product
\begin{equation*}
\langle u, v \rangle_{\J(\alpha,\beta)} := \int_{-1}^1 u(t) \, v(t) \, \J(\alpha,\beta)(t) \dd t
\end{equation*}
normalized according to \cite[Eq.~(4.1.1)]{szego:1975}
\begin{equation}\label{Jacobi1}
\PP_{n}^{(\alpha,\beta)}(1) = \binom{n+\alpha}{n}.
\end{equation}
Jacobi polynomials satisfy the Sturm--Liouville equation \cite[Eq.~(4.2.2)]{szego:1975}
\begin{equation*}
-\left( \J(\alpha+1,\beta+1) \, {\PP^{(\alpha,\beta)}_n}' \right)'
= \lambda^{(\alpha,\beta)}_n \, \J(\alpha,\beta) \, \PP^{(\alpha,\beta)}_n,
\end{equation*}
where
\begin{equation}\label{JacobiEigenvalue}
\lambda^{(\alpha,\beta)}_n := n (n+\alpha+\beta+1).
\end{equation}
By integration by parts we can express this equation in the well-known weak form:
\begin{equation}\label{JacobiWeakSL}
\left( \forall \, q \in \CC^1([-1,1]) \right) \quad
\left\langle {\PP^{(\alpha,\beta)}_n}', q' \right\rangle_{\J(\alpha+1,\beta+1)}
= \lambda^{(\alpha,\beta)}_n \, \left\langle \PP^{(\alpha,\beta)}_n, q \right\rangle_{\J(\alpha,\beta)}.
\end{equation}
Jacobi polynomials also satisfy the differentiation relation \cite[Eq.~(4.21.7)]{szego:1975}: For $n \geq 1$,
\begin{equation}\label{difshiftjac}
{\PP^{(\alpha,\beta)}_{n}}' = \frac{n + \alpha + \beta + 1}{2} \, \PP_{n-1}^{(\alpha+1,\beta+1)}.
\end{equation}
In what follows, we present four properties useful for the upcoming sections.

\begin{proposition}\label{parshift}
Let $\alpha,\beta>-1$ and $n\in\mathbb{N}_{0}$. Then,
\begin{equation*}
\norm{\PP^{(\alpha+1,\beta)}_{n}}_{\J(\alpha,\beta)}^{2} = \frac{2n+\alpha+\beta+2}{2(\alpha+1)}\norm{\PP_{n}^{(\alpha+1,\beta)}}_{\J(\alpha+1,\beta)}^{2}.
\end{equation*}
\begin{proof}
The polynomials $(1+\cdot) {P^{(\alpha+1,\beta)}_n}'$ and $n\,\PP^{(\alpha+1,\beta)}_n$ share the same leading term in their respective expansions with respect to the basis $\left\{(1+t)^k\right\}_{k=0}^{n}$, so their difference is a polynomial of degree less than or equal to $n-1$.
Therefore,
\begin{align*}
2n \norm{\PP^{(\alpha+1,\beta)}_n}_{\mathrm{J}(\alpha+1,\beta)}^2
&= 2\left\langle \PP^{(\alpha+1,\beta)}_n, n \, \PP^{(\alpha+1,\beta)}_n \right\rangle_{\mathrm{J}(\alpha+1,\beta)}\\
&= 2\left\langle \PP^{(\alpha+1,\beta)}_n, \J(0,1) {\PP^{(\alpha+1,\beta)}_n}' \right\rangle_{\mathrm{J}(\alpha+1,\beta)}\\
&= \int_{-1}^1 2 \PP^{(\alpha+1,\beta)}_n(t) {\PP^{(\alpha+1,\beta)}_n}'(t) (1-t)^{\alpha+1} (1+t)^{\beta+1} \, \mathrm{d}t\\
&= \int_{-1}^1 \left[ \PP^{(\alpha+1,\beta)}_n(t)^2 \right]' (1-t)^{\alpha+1} (1+t)^{\beta+1} \, \mathrm{d}t\\
&\stackrel{\text{IBP}}{=} -\int_{-1}^1 \PP^{(\alpha+1,\beta)}_n(t)^2 \left[ (1-t)^{\alpha+1} (1+t)^{\beta+1} \right]' \, \mathrm{d}t.
\end{align*}
Combining the above with the identity
\begin{equation*}
-\left[ \mathrm{J}(\alpha+1,\beta+1) \right]'
= 2(\alpha+1) \mathrm{J}(\alpha,\beta) - (\alpha+\beta+2) \mathrm{J}(\alpha+1,\beta),
\end{equation*}
we find that
\begin{multline*}
2n \norm{\PP^{(\alpha+1,\beta)}_n}_{\mathrm{J}(\alpha+1,\beta)}^2\\
= 2(\alpha+1) \norm{\PP^{(\alpha+1,\beta)}_n}_{\mathrm{J}(\alpha,\beta)}^2
- (\alpha+\beta+2) \norm{\PP^{(\alpha+1,\beta)}_n}_{\mathrm{J}(\alpha+1,\beta)}^2.
\end{multline*}
The equality above implies the desired result.
\end{proof}
\end{proposition}

\begin{proposition}\label{totalshift}
Let $\alpha,\beta>-1$ and $n\in\mathbb{N}_{0}$. Then,
\begin{equation*}
\norm{\PP^{(\alpha+1,\beta+1)}_{n}}_{\J(\alpha,\beta)}^{2} = \frac{(2n+\alpha+\beta+3)(\alpha+\beta+2)}{4(\alpha+1)(\beta+1)}\norm{\PP_{n}^{(\alpha+1,\beta+1)}}_{\J(\alpha+1,\beta+1)}^{2}.
\end{equation*}
\begin{proof}
We start by noticing that
\begin{multline}\label{JacobiWeightConnections}
-\big( (\beta-\alpha) + (\alpha+\beta+2) t \big) \, \mathrm{J}(\alpha+1,\beta+1)'(t)\\
= 4(\alpha+1)(\beta+1) \, \mathrm{J}(\alpha,\beta)(t) - (\alpha+\beta+2)^2 \, \mathrm{J}(\alpha+1,\beta+1)(t).
\end{multline}
Next, as the polynomials $n \, \PP^{(\alpha+1,\beta+1)}_n$ and
\begin{equation*}
\left(\frac{\beta - \alpha}{\alpha+\beta+2} + t\right) {\PP^{(\alpha+1,\beta+1)}_n}'(t)
\end{equation*}
share the same leading term in their respective expansions with respect to the basis
\begin{equation*}
\left\{1,\left(\frac{\beta - \alpha}{\alpha+\beta+2} + t\right),\dots,\left(\frac{\beta - \alpha}{\alpha+\beta+2} + t\right)^n \right\},
\end{equation*}
their difference is of degree less than or equal to $n-1$.
Therefore,
\begin{multline*}
2 n \norm{\PP^{(\alpha+1,\beta+1)}_n}_{\mathrm{J}(\alpha+1,\beta+1)}^2
= 2 \left\langle \PP^{(\alpha+1,\beta+1)}_n, n \, \PP^{(\alpha+1,\beta+1)}_n \right\rangle_{\mathrm{J}(\alpha+1,\beta+1)}\\
= 2 \int_{-1}^{1} \PP^{(\alpha+1,\beta+1)}_{n}(t) \left(\frac{\beta - \alpha}{\alpha+\beta+2} \! + \! t\right) {\PP^{(\alpha+1,\beta+1)}_n}'(t) \, \mathrm{J}(\alpha+1,\beta+1)(t)~\mathrm{d}t.
\end{multline*}
Using that
\begin{multline*}
\left[ \left(\frac{\beta - \alpha}{\alpha+\beta+2} \! + \! t\right) \PP^{(\alpha+1,\beta+1)}_n(t)^2 \right]'\\
= \left(\frac{\beta - \alpha}{\alpha+\beta+2} \! + \! t\right) 2 \PP^{(\alpha+1,\beta+1)}_n(t) {\PP^{(\alpha+1,\beta+1)}_n}'(t) + \PP^{(\alpha+1,\beta+1)}_n(t)^2.
\end{multline*}
It follows that
\begin{multline*}
(2 n + 1) \norm{\PP^{(\alpha+1,\beta+1)}_n}_{\mathrm{J}(\alpha+1,\beta+1)}^2\\
= \int_{-1}^1 \left[ \left(\frac{\beta - \alpha}{\alpha+\beta+2} \! + \! t\right) \PP^{(\alpha+1,\beta+1)}_n(t)^2 \right]' \mathrm{J}(\alpha+1,\beta+1)(t) \, \mathrm{d} t.
\end{multline*}
By integration by parts and using \eqref{JacobiWeightConnections} this results in
\begin{multline*}
(2 n + 1) \norm{\PP^{(\alpha+1,\beta+1)}_n}_{\mathrm{J}(\alpha+1,\beta+1)}^2\\
= \frac{4(\alpha+1)(\beta+1)}{\alpha+\beta+2} \norm{\PP^{(\alpha+1,\beta+1)}_n}_{\mathrm{J}(\alpha,\beta)}^2 - (\alpha+\beta+2) \norm{\PP^{(\alpha+1,\beta+1)}_n}_{\mathrm{J}(\alpha+1,\beta+1)}^2.
\end{multline*}
The equality above implies the desired result.
\end{proof}
\end{proposition}

\begin{proposition}\label{prop:double-shift}
Let $\alpha,\beta>-1$ and $n \in \mathbb{N}_{0}$.
Then,
\begin{equation*}
\norm{\PP^{(\alpha+2,\beta)}_{n}}_{\J(\alpha,\beta)}^{2} = \mathfrak{C}^{(\alpha,\beta)}(n)\norm{\PP_{n}^{(\alpha+2,\beta)}}_{\J(\alpha+2,\beta)}^{2}.
\end{equation*}
where $\mathfrak{C}^{(\alpha,\beta)}(n)$ is a polynomial in $n$ of degree $3$, given by
\begin{equation*}
\mathfrak{C}^{(\alpha,\beta)}(n) := \frac{\left(3 + \alpha + \beta + 2n\right) \left((3 + \alpha)(2 + \alpha + \beta) + 2(3 + \alpha + \beta)n + 2n^2\right)}{4(1 + \alpha)(2 + \alpha)(3 + \alpha)}.
\end{equation*}
\begin{proof}
As the polynomials $\left( 2n(n-1) + 4(\alpha+\beta+3) n \right) \PP^{(\alpha+2,\beta)}_n$ and
\begin{equation*}
2(1+t)^2 \, {\PP^{(\alpha+2,\beta)}_n}''(t) + 4(\alpha+\beta+3) (1+t) \, {\PP^{(\alpha+2,\beta)}_n}'(t)
\end{equation*}
share their principal terms, their difference is of degree lower than $n$, so
\begin{multline*}
\left( 2n(n-1) + 4(\alpha+\beta+3) n \right) \norm{\PP^{(\alpha+2,\beta)}_n}_{\J(\alpha+2,\beta)}^2\\
= \left\langle \PP^{(\alpha+2,\beta)}_n, 2 (1+t)^2 \, {\PP^{(\alpha+2,\beta)}_n}'' + 4(\alpha+\beta+3)(1+t) \, {\PP^{(\alpha+2,\beta)}_n}' \right\rangle_{\J(\alpha+2,\beta)}\\
= \int_{-1}^1 \biggl[ \Bigl( \bigl[ \PP^{(\alpha+2,\beta)}_n(t)^2 \bigr]'' - 2 \, {\PP^{(\alpha+2,\beta)}_n}'(t)^2 \Bigr) (1+t)^2\\ + 2(\alpha+\beta+3) \left[ \PP^{(\alpha+2,\beta)}_n(t)^2 \right]' (1+t) \biggr] \J(\alpha+2,\beta)(t) \dd t.
\end{multline*}
By integration by parts of the terms involving derivatives of the square of $\PP^{(\alpha+2,\beta)}_n$ (with vanishing boundary terms on account of the restriction $\alpha,\beta > -1$) and rearranging the term involving the square of the derivative of $\PP^{(\alpha+2,\beta)}_n$ with the help of the fact that $(1+t)^2 \, \J(\alpha+2,\beta)(t) = 2 \, \J(\alpha+2,\beta+1)(t) - \J(\alpha+3,\beta+1)(t)$, we can express the above as
\begin{multline}\label{double-shift-after-IBP}
2n(n+2\alpha+2\beta+5) \norm{\PP^{(\alpha+2,\beta)}_n}_{\J(\alpha+2,\beta)}^2\\
= \int_{-1}^1 \PP^{(\alpha+2,\beta)}_n(t)^2 \left[ \J(\alpha+2,\beta+2)''(t) - 2(\alpha+\beta+3) \J(\alpha+2,\beta+1)'(t) \right] \dd t\\
- 4 \norm{{\PP^{(\alpha+2,\beta)}_n}'}_{\J(\alpha+2,\beta+1)}^2 + 2 \norm{{\PP^{(\alpha+2,\beta)}_n}'}_{\J(\alpha+3,\beta+1)}^2.
\end{multline}
Let us note now that
\begin{multline}\label{longerJacobiWeightConnection}
\J(\alpha+2,\beta+2)'' - 2(\alpha+\beta+3) \, \J(\alpha+2,\beta+1)'\\
= 4(\alpha+1)(\alpha+2) \, \J(\alpha,\beta) - (\alpha+\beta+2)(\alpha+\beta+3) \, \J(\alpha+2,\beta).
\end{multline}
Further, by \eqref{difshiftjac}, ${\PP^{(\alpha+2,\beta)}_n}'$ is a scalar multiple of $\PP^{(\alpha+3,\beta+1)}_{n-1}$, so by \autoref{parshift},
\begin{equation}\label{double-shift-aux-1}
\norm{{\PP^{(\alpha+2,\beta)}_n}'}_{\J(\alpha+2,\beta+1)}^2
= \frac{2n+\alpha+\beta+3}{2(\alpha+3)} \norm{{\PP^{(\alpha+2,\beta)}_n}'}_{\J(\alpha+3,\beta+1)}^2
\end{equation}
and, by \eqref{JacobiWeakSL},
\begin{equation}\label{double-shift-aux-2}
\norm{{\PP^{(\alpha+2,\beta)}_n}'}_{\J(\alpha+3,\beta+1)}^2
= n(n+\alpha+\beta+3) \norm{\PP^{(\alpha+2,\beta)}_n}_{\J(\alpha+2,\beta)}^2.
\end{equation}
By combining \eqref{double-shift-after-IBP}, \eqref{longerJacobiWeightConnection}, \eqref{double-shift-aux-1} and \eqref{double-shift-aux-2} and rearranging, we obtain the desired result.
\end{proof}
\end{proposition}

\begin{proposition}\label{prop:Markov-sharpness-1d}
Let $\alpha,\beta>-1$ and $n\in\mathbb{N}_{0}$.
Then,
\begin{equation*}
\norm{{\PP_{n}^{(\alpha+1,\beta)}}'}_{\J(\alpha,\beta)}^{2} = \mathfrak{R}^{(\alpha,\beta)}(n)\norm{\PP_{n}^{(\alpha+1,\beta)}}_{\J(\alpha,\beta)}^{2},
\end{equation*}
where $\mathfrak{R}^{(\alpha,\beta)}(n)$ is the polynomial of degree $4$ in $n$ given by
\begin{equation*}
\mathfrak{R}^{(\alpha,\beta)}(n) = \frac{2(\beta+1)\left(\lambda_{n}^{(\alpha+1,\beta)}\right)^2 + (\alpha+1)(\alpha + \beta + 3)(\alpha+\beta+4) \lambda_{n}^{(\alpha+1,\beta)}}{4(\alpha+2)(\alpha+3)(\beta+1)}.
\end{equation*}
\begin{proof}
If $n=0$, both sides vanish, so we assume $n\geq 1$.
Using that $\J(\alpha,\beta) = \frac{1}{2}\left[\J(\alpha+1,\beta)+\J(\alpha,\beta+1)\right]$ we have the identity
\begin{equation}\label{eq:split}
\norm{\PP_{n-1}^{(\alpha+2,\beta+1)}}_{\J(\alpha,\beta)}^{2} = \frac{1}{2}\left(\norm{\PP_{n-1}^{(\alpha+2,\beta+1)}}_{\J(\alpha+1,\beta)}^{2} + \norm{\PP_{n-1}^{(\alpha+2,\beta+1)}}_{\J(\alpha,\beta+1)}^{2}\right).
\end{equation}
For the first term, \autoref{totalshift} with parameters $(\alpha+1,\beta)$ and degree $n-1$ gives
\begin{equation*}
\norm{\PP_{n-1}^{(\alpha+2,\beta+1)}}_{\J(\alpha+1,\beta)}^{2} = \frac{(2n+\alpha+\beta+2)(\alpha+\beta+3)}{4(\alpha+2)(\beta+1)}\norm{\PP_{n-1}^{(\alpha+2,\beta+1)}}_{\J(\alpha+2,\beta+1)}^{2}.
\end{equation*}
For the second term, \autoref{prop:double-shift} with parameters $(\alpha,\beta+1)$ and degree $n-1$ gives
\begin{equation*}
\norm{\PP_{n-1}^{(\alpha+2,\beta+1)}}_{\J(\alpha,\beta+1)}^{2} = \mathfrak{C}^{(\alpha,\beta+1)}(n-1)\norm{\PP_{n-1}^{(\alpha+2,\beta+1)}}_{\J(\alpha+2,\beta+1)}^{2}.
\end{equation*}
The factor $\mathfrak{C}^{(\alpha,\beta+1)}(n-1)$ can we written as
\begin{equation*}
\mathfrak{C}^{(\alpha,\beta+1)}(n-1) = \frac{(2n+\alpha+\beta+2)\left(2\lambda_{n}^{(\alpha+1,\beta)} + (\alpha+1)(\alpha+\beta+3)\right)}{4(\alpha+1)(\alpha+2)(\alpha+3)}.
\end{equation*}
Substituting these two identities into \eqref{eq:split} and using the fact that
\begin{equation*}
{\PP_{n}^{(\alpha+1,\beta)}}' \stackrel{\eqref{difshiftjac}}{=} \frac{n+\alpha+\beta+2}{2} \, \PP_{n-1}^{(\alpha+2,\beta+1)},
\end{equation*}
we obtain
\begin{multline}\label{eq:num}
\norm{{\PP_{n}^{(\alpha+1,\beta)}}'}_{\J(\alpha,\beta)}^{2} = \frac{2n+\alpha+\beta+2}{8(\alpha+2)}\\
\times \biggl(\frac{\alpha+\beta+3}{\beta+1}
 + \frac{2\lambda_{n}^{(\alpha+1,\beta)} + (\alpha+1)(\alpha+\beta+3)}{(\alpha+1)(\alpha+3)}\biggr) \norm{{\PP_{n}^{(\alpha+1,\beta)}}'}_{\J(\alpha+2,\beta+1)}^{2}.
\end{multline}
Next, using \eqref{JacobiWeakSL} with $\alpha+1$ instead of $\alpha$ and $q = \PP_{n}^{(\alpha+1,\beta)}$ yields 
\begin{equation*}
\norm{{\PP_{n}^{(\alpha+1,\beta)}}'}_{\J(\alpha+2,\beta+1)}^{2} = \lambda_{n}^{(\alpha+1,\beta)} \norm{\PP_{n}^{(\alpha+1,\beta)}}_{\J(\alpha+1,\beta)}^{2},
\end{equation*}
and \autoref{parshift} provides
\begin{equation*}
\norm{\PP_{n}^{(\alpha+1,\beta)}}_{\J(\alpha+1,\beta)}^{2} = \frac{2(\alpha+1)}{2n+\alpha+\beta+2}\norm{\PP_{n}^{(\alpha+1,\beta)}}_{\J(\alpha,\beta)}^{2}.
\end{equation*}
Combining the last two identities with \eqref{eq:num} and rearranging we deduce the desired result.
\end{proof}
\end{proposition}

\section{Orthogonal Polynomial Spaces}\label{sec:Orthogonalpolyspace}
\subsection{Orthogonal polynomials with respect to a Jacobi-type weight}
Let $\mathcal{V}_{n}^{\Gam}$ be the space of the orthogonal polynomials of degree $n$ with respect to the weight $\mathrm{W}_{\Gam}$; i.e.,
\begin{equation}\label{SpaceOP}
\mathcal{V}_{n}^{\Gam} := \{ p\in\Pi_{n}^{d}| ~ (\forall q\in\Pi_{n-1}^{d})\ \langle p,q\rangle_{\Gam} = 0\}.
\end{equation}
For $n<0$, from the convention $\Pi_{n}^{d} = \{0\}$ it follows that $\mathcal{V}_{n}^{\Gam} = \{0\}$.
There holds (cf.\ \cite[Sec.~3.1]{DunklXu:2014})
\begin{equation}\label{eq:dimVn}
\operatorname{dim}(\mathcal{V}_{n}^{\Gam}) = \binom{n+d-1}{n}.
\end{equation}
Let $S_{n}^{\Gam}$ and $\proj_{n}^{\Gam}$ denote the orthogonal projections from $\LL_{\Gam}^{2}$ onto $\Pi_{n}^{d}$ and from $\LL_{\Gam}^{2}$ onto $\mathcal{V}_{n}^{\Gam}$, respectively.
From \cite[Th.~3.2.18]{DunklXu:2014},
\begin{equation}\label{eq:orthogonal-decomposition-of-polynomials}
\Pi_{n}^{d} = \bigoplus_{k=0}^{n}\mathcal{V}_{k}^{\Gam}\quad \text{and} \quad \mathrm{L}_{\Gam}^{2} = \bigoplus_{k=0}^{\infty}\mathcal{V}_{k}^{\Gam},
\end{equation}
whence
\begin{equation}\label{projectionform}
(\forall n\in\mathbb{N}_{0})\quad S_{n}^{\Gam} = \sum_{k=0}^{n}\proj_{k}^{\Gam}\quad \text{and}\quad (\forall u\in\LL_{\Gam}^{2})\quad u = \sum_{k=0}^{\infty}\proj_{k}^{\Gam}(u).
\end{equation}
In what follows we will make heavy use of the definition \eqref{weight} so that, e.g.,
\begin{equation}\label{eq:Wej}
\W_{e_{j}+e_{d+1}}(x) = x_{j}(1-\abs{x}).
\end{equation}
We introduce the second-order differential operator
\begin{equation}\label{L}
\mathcal{L}^{\Gam}q
:= -\frac{1}{\W_{\Gam}} \left[ \sum_{i=1}^d \partial_i\left( \W_{\Gam+e_i+e_{d+1}} \partial_i q \right) + \sum_{1\leq i<j \leq d} \partial_{i,j} \left( \W_{\Gam+e_i+e_j} \partial_{i,j} q \right) \right],
\end{equation}
where $\partial_{i,j} = \partial_j - \partial_i$.
For $\Gam\in(-1,\infty)^{d+1}$, members of the space of orthogonal polynomials $\mathcal{V}_{n}^{\Gam}$ satisfy the second-order Sturm--Liouville problem (cf.\ \cite[Sec.~5.3]{DunklXu:2014}, \cite[Eq.~(2.1)]{GeXu:2024}).
\begin{equation}\label{SL-O-strong}
\left(\forall p_{n}\in\mathcal{V}_{n}^{\Gam}\right)\quad \mathcal{L}^{\Gam}(p_{n}) = \lambda_{n}^{\Gam} \, p_{n},
\end{equation}
where, $\lambda_{n}^{\Gam} = n(n+\abs{\Gam}+d)$.

In what follows, we provide a characterization of the space of orthogonal polynomials $\mathcal{V}_{n}^{\Gam}$ on the $d$-simplex in terms of mapped tensor products of spaces of orthogonal polynomials on the $(d-1)$-simplex and orthogonal polynomials on the interval $(-1,1)$. To this end, we introduce the nonlinear map $\Psi_{d}:\tri^{d-1}\times(-1,1)\longrightarrow \tri^{d}$ defined by
\begin{equation}\label{Psi-d}
\Psi_{d}(\bd{\zeta},\eta) := \left( \frac{1-\eta}{2}\bd{\zeta},\frac{(1-\eta)(1-\abs{\bd{\zeta}})}{2}\right),
\end{equation}
whose Jacobian determinant is
\begin{equation*}
\det\left[\operatorname{D}(\Psi_{d})(\bd{\zeta},\eta)\right] = 2^{-d}(1-\eta)^{d-1}.
\end{equation*}
\begin{definition}\label{def:tensor-mathbb-T}
Let $d\in\mathbb{N}_{\geq2}$, $l\in\mathbb{N}_{0}$, $k\in[l]_{0}$ and $\Gam = (\tilde{\Gam},\gamma_{d+1})\in(-1,\infty)^{d+1}$.
We define the linear operator $\mathbb{T}_{l,k}^{\Gam}:\Pi_{k}^{d-1}\longrightarrow \Pi_{l}^{d}$ as
\begin{equation*}
\mathbb{T}_{l,k}^{\Gam}(q)(x) := \abs{x}^{k}q\left(\frac{x_{1}}{\abs{x}},\dots,\frac{x_{d-1}}{\abs{x}}\right)\PP_{l-k}^{(2k+\abs{\tilde{\Gam}}+d-1,\gamma_{d+1})}(1-2\abs{x}).
\end{equation*}
\end{definition}
\begin{proposition}\label{prop:poly-tensor}
Let $d\in\mathbb{N}_{\geq3}$, $n,m\in\mathbb{N}_{0}$, $k\in[n]_{0}$, $l\in[m]_{0}$, $\Gam = (\tilde{\Gam},\gamma_{d+1})\in(-1,\infty)^{d+1}$ and $\boldsymbol{\kappa} = (\tilde{\boldsymbol{\kappa}},\kappa_{d+1})\in(-1,\infty)^{d+1}$.
Then, for all $p\in\Pi_{k}^{d-1}$ and $q\in\Pi_{l}^{d-1}$,
\begin{multline*}
\langle \mathbb{T}_{n,k}^{\Gam}(p),\mathbb{T}_{m,l}^{\Gam}(q) \rangle_{\boldsymbol{\kappa}} = 2^{-\abs{\boldsymbol{\kappa}}-k-l-d}\langle p,q\rangle_{\tilde{\boldsymbol{\kappa}}}\\
\times \left\langle \PP_{n-k}^{(2k+\abs{\tilde{\Gam}}+d-1,\gamma_{d+1})},\PP_{m-l}^{(2l+\abs{\tilde{\Gam}}+d-1,\gamma_{d+1})} \right\rangle_{\J(k+l+\abs{\tilde{\boldsymbol{\kappa}}}+d-1,\kappa_{d+1})}.
\end{multline*}
\begin{proof}
The result follows directly by applying the change of variables $(x_{1},\dots,x_{d}) = \Psi_{d}(\bd{\zeta},\eta)$ to the left-hand side of the identity.
\end{proof}
\begin{remark}\label{ortho-mathbb-T}
From \autoref{prop:poly-tensor} and the orthogonality properties of Jacobi polynomials, we deduce that if $p\in\mathcal{V}_{k}^{\tilde{\Gam}}$ and $q\in\mathcal{V}_{l}^{\tilde{\Gam}}$, then $\mathbb{T}^{\Gam}_{n,k}(p)\perp_{\W_{\Gam}}\mathbb{T}^{\Gam}_{m,l}(q)$ when $(n,k)\neq(m,l)$.
\end{remark}
\begin{remark}\label{rem:injectivity-tensorization-mathbb-T}
In view of \autoref{prop:poly-tensor} and the fact that the involved Jacobi polynomial is not null (cf.\ \eqref{Jacobi1}), it follows that the operator $\mathbb{T}_{l,k}^{\Gam}$ introduced in \autoref{def:tensor-mathbb-T} is injective.
\end{remark}
\end{proposition}
\begin{proposition}\label{prop:decomposition-V}
Let $d\in\mathbb{N}_{\geq3}$, $n\in\mathbb{N}_{0}$, and $\Gam =(\tilde{\Gam},\gamma_{d+1})\in (-1,\infty)^{d+1}$.
Then,
\begin{equation}\label{eq:decom-V-space}
\mathcal{V}_{n}^{\Gam} = \bigoplus_{k=0}^{n}\mathbb{T}_{n,k}^{\Gam}\left(\mathcal{V}_{k}^{\tilde{\Gam}}\right).
\end{equation}
Moreover, this decomposition is $\W_{\Gam}$-orthogonal.
\begin{proof}
Let $i,j\in[n]_{0}$ be such that $i\neq j$.
By \autoref{ortho-mathbb-T}, we have that
\begin{equation*}
\mathbb{T}^{\Gam}_{n,i}(\mathcal{V}_{i}^{\tilde{\Gam}})\perp_{\W_{\Gam}}\mathbb{T}^{\Gam}_{n,j}(\mathcal{V}_{j}^{\tilde{\Gam}}).
\end{equation*}
Consequently, the decomposition on the right-hand side of \eqref{eq:decom-V-space} is $\W_{\Gam}$-orthogonal.
This fact, together with \autoref{rem:injectivity-tensorization-mathbb-T} and \eqref{eq:dimVn}, leads to
\begin{equation}\label{eq:dim-oplus-mathbb-T}
\operatorname{dim}\mathopen{}\left(\bigoplus_{k=0}^{n}\mathbb{T}^{\Gam}_{n,k}(\mathcal{V}^{\tilde{\Gam}}_{k})\mathclose{}\right)
= \sum_{k=0}^{n}\operatorname{dim}\mathopen{}\left(\mathbb{T}_{n,k}^{\Gam}(\mathcal{V}_{k}^{\tilde{\Gam}})\mathclose{}\right)
= \sum_{k=0}^{n}\binom{k+d-2}{k} = \binom{n+d-1}{n}.
\end{equation}
In view of \eqref{eq:dimVn} and the above, it follows that $\operatorname{dim}\mathopen{}\left(\bigoplus_{k=0}^{n}\mathbb{T}^{\Gam}_{n,k}(\mathcal{V}^{\tilde{\Gam}}_{k})\mathclose{}\right) = \operatorname{dim}(\mathcal{V}_{n}^{\Gam})$.
We shall now show that
\begin{equation}\label{eq:incluoplusT}
\bigoplus_{k=0}^{n}\mathbb{T}^{\Gam}_{n,k}(\mathcal{V}^{\tilde{\Gam}}_{k})\subset \mathcal{V}_{n}^{\Gam},
\end{equation}
which, combined with \eqref{eq:dim-oplus-mathbb-T}, establishes the equality \eqref{eq:decom-V-space}.
For $m\in\mathbb{N}_{0}$, we first observe that \autoref{ortho-mathbb-T} and \eqref{eq:dim-oplus-mathbb-T} imply
\begin{equation*}
\dim\mathopen{}\left(\bigoplus_{l=0}^{m}\bigoplus_{r=0}^{l}\mathbb{T}^{\Gam}_{l,r}\left(\mathcal{V}_{r}^{\tilde{\Gam}}\right)\mathclose{}\right)
= \sum_{l=0}^{m}\binom{l+d-1}{l} = \binom{m+d}{m} = \operatorname{dim}(\Pi_{m}^{d}).
\end{equation*}
Since $\mathbb{T}^{\Gam}_{l,r}\left(\Pi_{r}^{d-1}\right)\subset\Pi_{l}^{d}$ (cf.~\autoref{def:tensor-mathbb-T}), we deduce that
\begin{equation}\label{eq:descomPoly}
\Pi_{m}^{d} = \bigoplus_{l=0}^{m}\bigoplus_{r=0}^{l}\mathbb{T}^{\Gam}_{l,r}\left(\mathcal{V}_{r}^{\tilde{\Gam}}\right).
\end{equation}
Thus, given a polynomial $p_{n} \in \bigoplus_{k=0}^{n}\mathbb{T}_{n,k}^{\Gam}(\mathcal{V}_{k}^{\tilde{\Gam}}) \subset \Pi_{n}^{d}$ and any polynomial $q \in \Pi_{n-1}^{d}$, we can expand the latter according to \eqref{eq:descomPoly} and so we can apply \autoref{ortho-mathbb-T} to conclude that $\langle p_{n}, q_{n-1} \rangle_{\Gam} = 0$.
\end{proof}
\end{proposition}
\subsection{Orthogonal polynomials with respect to Jacobi-type weight with a ridge singularity}\label{subsec:Dep-OP}
Given $d \geq 2$, the set $\{ x \in \overline{\tri^d} \mid x_1 = x_2 = 0 \}$ is a \emph{ridge} or $(d-2)$-face of the closed simplex $\overline{\tri^d}$; in particular, in the $d = 2$ case it is a vertex and in the $d = 3$ case it is an edge.
Given $d \geq 2$ and $\Gam \in \mathbb{R}^{d+1}$ we define a Jacobi-type weight function $\W_{\mathrm{ridge},\Gam} \colon \tri^{d}\longrightarrow\mathbb{R}$ involving an additional singularity near said ridge by
\begin{equation}\label{eq:depleted-weight}
\left(\forall\,x\in\tri^{d}\right)\quad \W_{\mathrm{ridge},\Gam}(x) := \frac{\W_{\Gam}(x)}{x_{1}+x_{2}}.
\end{equation}
For those $\Gam\in(-1,\infty)^{d+1}$ satisfying $\W_{\mathrm{ridge},\Gam}\in\LL^{1}(\tri^{d})$ (sufficient conditions are given in \autoref{NewWeight} below), we define the space
\begin{equation*}
\LL^{2}(\W_{\mathrm{ridge},\Gam}) := \{\W_{\mathrm{ridge},\Gam}^{-1/2}f\mid f\in\LL^{2}(\tri^{d})\},
\end{equation*}
whose natural inner product is
\begin{equation*}
\left(\forall f,g\in \LL^{2}(\W_{\mathrm{ridge},\Gam})\right) \quad \langle f,g\rangle_{\mathrm{ridge},\Gam} := \int_{\tri^{d}}f(x) g(x) \W_{\mathrm{ridge},\Gam}(x) \, \textup{d}x.
\end{equation*}
We denote by $\norm{\cdot}_{\mathrm{ridge},\Gam}$ the norm induced by the inner product $\langle\cdot,\cdot\rangle_{\mathrm{ridge},\Gam}$.
The condition $\W_{\mathrm{ridge},\Gam}\in\LL^{1}(\tri^{d})$ ensures that every $d$-variate polynomial belongs to $\LL^{2}(\W_{\mathrm{ridge},\Gam})$.
Given $n\in\mathbb{N}_{0}$, we denote by $\mathcal{V}_{n}^{\mathrm{ridge},\Gam}$ the space of $\LL^{2}(\W_{\mathrm{ridge},\Gam})$-orthogonal polynomials of degree $n$; that is,
\begin{equation}\label{DepletedOP}
\mathcal{V}_{n}^{\mathrm{ridge},\Gam} := \left\{p\in\Pi_{n}^{d}\mid \left(\forall q\in\Pi_{n-1}^{d}\right)~ \langle p,q\rangle_{\mathrm{ridge},\Gam} = 0\right\}.
\end{equation}
We adopt the convention $\mathcal{V}_{n}^{\mathrm{ridge},\Gam} = \{0\}$ for $n<0$.
There holds (cf.~\cite[Sec.~3.1]{DunklXu:2014})
\begin{equation}\label{dimOPS}
\left(\forall n\in\mathbb{N}_{0}\right)\quad \operatorname{dim}\left(\mathcal{V}_{n}^{\mathrm{ridge},\Gam}\right) = \binom{n+d-1}{n}.
\end{equation}
Let $\proj_{n}^{\mathrm{ridge},\Gam}$ denote the orthogonal projection from $\LL^{2}(\W_{\mathrm{ridge},\Gam})$ onto $\mathcal{V}_{n}^{\mathrm{ridge},\Gam}$.
From \cite[Th.~3.2.18]{DunklXu:2014},
\begin{equation}\label{eq:dep-orthogonal-decomposition-of-polynomials}
\Pi_{n}^{d} = \bigoplus_{k=0}^{n}\mathcal{V}_{k}^{\mathrm{ridge},\Gam}\quad \textup{and}\quad \LL^{2}(\W_{\mathrm{ridge},\Gam}) = \bigoplus_{k=0}^{\infty}\mathcal{V}_{k}^{\mathrm{ridge},\Gam},
\end{equation}
whence
\begin{equation}\label{projectionformdep}
(\forall u\in\LL^{2}(\W_{\mathrm{ridge},\Gam}))\quad u = \sum_{k=0}^{\infty}\proj_{k}^{\mathrm{ridge},\Gam}(u).
\end{equation}
In the rest of this section provide a characterization of the space of orthogonal polynomials $\mathcal{V}_{n}^{\mathrm{ridge},\Gam}$ analogous to the characterization of $\mathcal{V}_n^{\Gam}$ given in \autoref{prop:decomposition-V}.
Such characterization will later enable an induction argument on the dimension to prove \autoref{lem:slice-depleted-Schur}, which is an immediate precursor of \eqref{crude-depleted-Schur}.

We shall again make use of the nonlinear map $\Psi_{d}$ of \eqref{Psi-d} but we shall also require, in the case of the triangle, the map $\Psi_{2} \colon (-1,1)\times(-1,1) \to \tri^{2}$ defined as
\begin{equation}\label{Psi}
\Psi_{2}(\zeta,\eta) := \left( \frac{(1-\eta)(1-\zeta)}{4},\frac{(1-\eta)(1+\zeta)}{4} \right),
\end{equation}
whose Jacobian determinant is
\begin{equation*}
\det\left[\operatorname{D}(\Psi_{2})(\zeta,\eta)\right] = 2^{-3}(1-\eta).
\end{equation*}
Let us further define the special parameter set
\begin{equation*}
\Lambda^{d} := \{\Kap\in\mathbb{R}^{d+1}\mid \Kap\in(-1,\infty)^{d+1}\ \wedge \ \kappa_{1}+\kappa_{2}>-1\}.
\end{equation*}

\begin{proposition}\label{NewWeight}
Let $d \geq 2$ and let $\Gam\in\Lambda^d$.
Then, $\W_{\mathrm{ridge},\Gam}\in\LL^{1}(\tri^{d})$.
\begin{proof}
We prove the result by induction on the dimension $d$.

For $d=2$, let $\Gam\in\Lambda^{2}$.
Using the change of variables $(x_{1},x_{2}) = \Psi_{2}(\zeta,\eta)$ and known integrals (cf.\ \cite[Eqs.~(1.7.4--5)]{szego:1975}) we have
\begin{multline*}
\int_{\tri^{2}}\W_{\mathrm{ridge},\Gam}(x)\ \textup{d}x = 2^{-(2\gamma_{1}+2\gamma_{2}+\gamma_{3}+2)}\int_{-1}^{1} \J(\gamma_{1}+\gamma_{2},\gamma_{3})(\eta)\ \textup{d}\eta \int_{-1}^{1} \J(\gamma_{1},\gamma_{2})(\zeta)\ \textup{d}\zeta\\
= \frac{\Gamma(\gamma_{1}+1)\Gamma(\gamma_{2}+1)\Gamma(\gamma_{3}+1)}{(\gamma_{1}+\gamma_{2}+1) \Gamma(\gamma_{1}+\gamma_{2}+\gamma_{3}+2)}.
\end{multline*}
Since $\Gam\in\Lambda^{2}$, we have $\gamma_{1},\gamma_{2},\gamma_{3}>-1$ and $\gamma_{1}+\gamma_{2}>-1$, which ensures that the expression above is finite.
Thus, $\W_{\mathrm{ridge},\Gam}\in\LL^{1}(\tri^{2})$.

Assume for the induction hypothesis that for some $d\in\mathbb{N}_{\geq 2}$ the statement holds for any parameter vector in $\Lambda^{d}$.
Now, let $\Gam := (\gamma_{1},\dots,\gamma_{d+1},\gamma_{d+2})\in\Lambda^{d+1}$.
Let $\tilde{\Gam} := (\gamma_{1},\dots,\gamma_{d+1})$, whence $\Gam = (\tilde{\Gam},\gamma_{d+2})$.
Performing the change of variables $(x_{1},\dots,x_{d+1}) = \Psi_{d+1}(\bd{\zeta},\eta)$, we obtain
\begin{equation*}
\int_{\tri^{d+1}}\W_{\mathrm{ridge},\Gam}(x)\ \textup{d}x = \frac{\Gamma(|\tilde{\Gam}|+d)\Gamma(\gamma_{d+2}+1)}{\Gamma(\abs{\Gam}+d+1)}\int_{\tri^{d}}\W_{\mathrm{ridge},\tilde{\Gam}}(x)\ \textup{d}x.
\end{equation*}
The condition $\Gam\in\Lambda^{d+1}$ implies that $\tilde{\Gam}\in\Lambda^{d}$, so, by our induction hypothesis, the integral on the right hand side above is finite.
Consequently, the integral on the left-hand side is also finite, which proves that $\W_{\mathrm{ridge},\Gam}\in\LL^{1}(\tri^{d+1})$.
This completes the induction.
\end{proof}
\end{proposition}
We now provide a basis for the space of $\W_{\mathrm{ridge},\Gam}$-orthogonal polynomials on the triangle $\tri^{2}$, which is an affine mapping away from a particular case of the basis of the $P^{(a,b,c,d)}_{n,k}$ of \cite{OlverTownsendVasil:2020} and is also closely connected with known bases of $\W_{\Gam}$-orthogonal polynomials \cite[Sec.~2.4]{DunklXu:2014}.
\begin{definition}\label{basisdepletedW}
Let $\Gam\in\Lambda^{2}$, $n\in\mathbb{N}_{0}$ and $j\in[n]_{0}$.
We define the family of polynomials $\{ \PP_{n,j}^{\mathrm{ridge},\Gam} \}_{j=0}^{n}$ by
\begin{equation*}
\PP_{n,j}^{\mathrm{ridge},\Gam}(x_{1},x_{2}) := (x_{1}+x_{2})^{j} \, \PP_{j}^{(\gamma_{1},\gamma_{2})}\left(\frac{x_{2}-x_{1}}{x_{1}+x_{2}}\right)\PP^{(2j+\gamma_{1}+\gamma_{2},\gamma_{3})}_{n-j}(1-2x_{1}-2x_{2}).
\end{equation*}
\end{definition}

It is readily apparent that, for all $\Gam \in \Lambda^{2}$, $n \in \mathbb{N}_0$ and $j \in [n]_0$,
\begin{equation}\label{depletedBasisMemberDegree}
\PP_{n,j}^{\mathrm{ridge},\Gam}\in\Pi^{2}_{n}\backslash\Pi_{n-1}^{2}.
\end{equation}

\begin{proposition}\label{orthobreve}
Let $\Gam\in \Lambda^{2}$ and $n\in\mathbb{N}_{0}$.
Then, $\{ \PP_{n,k}^{\mathrm{ridge},\Gam} \}_{k=0}^n$ is a $\W_{\mathrm{ridge},\Gam}$-orthogonal basis of $\mathcal{V}_{n}^{\mathrm{ridge},\Gam}$.
\begin{proof}
Let $n,m\in\mathbb{N}_{0}$, $l\in[n]_{0}$ and $j\in[m]_{0}$.
Performing the change of variables $(x_{1},x_{2}) = \Psi_{2}(\zeta,\eta)$ we have
\begin{multline}\label{eq:decom-depled-weight-2d}
\langle \PP_{m,j}^{\mathrm{ridge},\Gam},\PP_{n,l}^{\mathrm{ridge},\Gam} \rangle_{\mathrm{ridge},\Gam} = 2^{-(2+j+l+2\gamma_{1}+2\gamma_{2}+\gamma_{3})}
\left\langle \PP_{j}^{(\gamma_{1},\gamma_{2})},\PP_{l}^{(\gamma_{1},\gamma_{2})} \right\rangle_{\J(\gamma_{1},\gamma_{2})}\\
\times\left\langle \PP_{m-j}^{(2j+\gamma_{1}+\gamma_{2},\gamma_{3})},\PP_{n-l}^{(2l+\gamma_{1}+\gamma_{2},\gamma_{3})} \right\rangle_{\J(j+l+\gamma_{1}+\gamma_{2},\gamma_{3})}.
\end{multline}
Given $N \in \mathbb{N}_0$, from \eqref{depletedBasisMemberDegree}, \eqref{eq:decom-depled-weight-2d} and the orthogonality properties of Jacobi polynomials we find that $A_N := \{ \PP_{n,k}^{\mathrm{ridge},\Gam} \mid n \in [N]_0,\ k \in [n]_0 \}$ is a mutually orthogonal subset of $\Pi_{N}^{2}$; as $\# A_N = \dim(\Pi_{N}^{2})$, $\vspan(A_N) = \Pi_{N}^{2}$.
Given $n \in \mathbb{N}_0$ it follows from \eqref{depletedBasisMemberDegree} and \eqref{eq:decom-depled-weight-2d} and the fact that $\Pi_{n-1}^{2} = \vspan(A_{n-1})$, that the mutually $\W_{\mathrm{ridge},\Gam}$-orthogonal set $B_n := \{\PP_{n,k}^{\mathrm{ridge},\Gam}\}_{k=0}^n$ is contained in $\Pi_{n}^{2}\setminus \Pi_{n-1}^2$ and is $\W_{\mathrm{ridge},\Gam}$-orthogonal to $\Pi_{n-1}^{2}$, whence $\vspan(B_n) \subset \mathcal{V}_{n}^{\mathrm{ridge},\Gam}$; as $\# B_n = \dim(\mathcal{V}_{n}^{\mathrm{ridge},\Gam})$ (cf.\ \eqref{dimOPS}), $\vspan(B_n) = \mathcal{V}_{n}^{\mathrm{ridge},\Gam}$.
\end{proof}
\end{proposition}

We now show that the $\W_{\mathrm{ridge},\Gam}$-orthogonal basis of $\mathcal{V}_{n}^{\mathrm{ridge},\Gam+e_{1}+e_{2}}$ given by \autoref{orthobreve} is also $\W_{\Gam}$-orthogonal.

\begin{proposition}\label{orthw}
Let $\Gam\in(-1,\infty)^{3}$ and $n\in\mathbb{N}_{0}$.
Then, the set $\{ \PP_{n,k}^{\mathrm{ridge},\Gam+e_{1}+e_{2}} \}_{k=0}^n$ is mutually $\W_{\Gam}$-orthogonal.
\begin{proof}
Let $n,m\in\mathbb{N}_{0}$, $l\in[n]_{0}$ and $j\in[m]_{0}$.
Performing the change of variables $(x_{1},x_{2}) = \Psi_{2}(\zeta,\eta)$ we have
\begin{multline}\label{eq:decom-inner-product-2d}
\left\langle \PP_{m,j}^{\mathrm{ridge},\Gam+e_{1}+e_{2}},\PP_{n,l}^{\mathrm{ridge},\Gam+e_{1}+e_{2}} \right\rangle_{\Gam} = 2^{-(3+j+l+2\gamma_{1}+2\gamma_{2}+\gamma_{3})}\\
\times\left\langle \PP_{j}^{(\gamma_{1}+1,\gamma_{2}+1)},\PP_{l}^{(\gamma_{1}+1,\gamma_{2}+1)} \right\rangle_{\J(\gamma_{1},\gamma_{2})}\\
\times\left\langle \PP_{m-j}^{(2j+\gamma_{1}+\gamma_{2}+2,\gamma_{3})},\PP_{n-l}^{(2l+\gamma_{1}+\gamma_{2}+2,\gamma_{3})} \right\rangle_{\J(j+l+\gamma_{1}+\gamma_{2}+1,\gamma_{3})}.
\end{multline}
Supose that $n=m$ and $j < l$; that is, $j+1 \leq l$.
Then, the second inner product on the right-hand side of \eqref{eq:decom-inner-product-2d} can be rewritten as
\begin{multline*}
\left\langle \PP_{n-j}^{(2j+\gamma_{1}+\gamma_{2}+2,\gamma_{3})},\PP_{n-l}^{(2l+\gamma_{1}+\gamma_{2}+2,\gamma_{3})} \right\rangle_{\J(j+l+\gamma_{1}+\gamma_{2}+1,\gamma_{3})}\\
= \left\langle \PP_{n-j}^{(2j+\gamma_{1}+\gamma_{2}+2,\gamma_{3})},\J(l-j-1,0)\PP_{n-l}^{(2l+\gamma_{1}+\gamma_{2}+2,\gamma_{3})}  \right\rangle_{\J(2j+\gamma_{1}+\gamma_{2}+2,\gamma_{3})}.
\end{multline*}
Since $l - j - 1\geq 0$, the function $\J(l-j-1,0)\PP_{n-l}^{(2l+\gamma_{1}+\gamma_{2}+2,\gamma_{3})}$ is a polynomial of degree $n-j-1$.
Hence, by the orthogonality of the Jacobi polynomial, the above inner product vanishes.
The case in which $l < j$ is completely analogous.
In summary, when $n = m$ and $j \neq l$, the inner product \eqref{eq:decom-inner-product-2d} vanishes.
\end{proof}
\end{proposition}

Now we define a family of operators $\mathcal{T}_{l,j}^{\Gam}$, analogous to the $\mathbb{T}_{l,k}^{\Gam}$ of \autoref{def:tensor-mathbb-T}, whose purpose is turning $\W_{\mathrm{ridge},\Gam}$-orthogonal polynomials into analogues on higher dimension.

\begin{definition}\label{def:tensorization}
Let $d\in\mathbb{N}_{\geq 3}$, $l\in\mathbb{N}_{0}$, $k\in[l]_{0}$, and $\Gam := (\tilde{\Gam},\gamma_{d+1})\in \Lambda^{d}$.
We define the linear operator $\mathcal{T}^{\Gam}_{l,k}:\ \Pi_{k}^{d-1}\longrightarrow\Pi_{l}^{d}$ as
\begin{equation*}
\mathcal{T}^{\Gam}_{l,k}(q)(x) := \abs{x}^{k}q\left( \frac{x_{1}}{\abs{x}},\dots,\frac{x_{d-1}}{\abs{x}} \right)\PP_{l-k}^{(2k+\abs{\tilde{\Gam}}+d-2,\gamma_{d+1})}(1-2\abs{x}).
\end{equation*}
\end{definition}

\begin{proposition}\label{prop:normT}
Let $d\in\mathbb{N}_{\geq 3}$, $n,m\in\mathbb{N}_{0}$, $k\in[n]_{0}$, $l\in[m]_{0}$, and $\Gam := (\tilde{\Gam},\gamma_{d+1}) \in \Lambda^{d}$.
Then,
\begin{multline*}
\langle \mathcal{T}^{\Gam}_{n,k}(p),\mathcal{T}^{\Gam}_{m,l}(q) \rangle_{\mathrm{ridge},\Gam} = 2^{-(k+l+\abs{\Gam}+d-1)}\langle p,q \rangle_{\mathrm{ridge},\tilde{\Gam}}\\
\times\left\langle \PP_{n-k}^{(2k+\abs{\tilde{\Gam}}+d-2,\gamma_{d+1})},\PP_{m-l}^{(2l+\abs{\tilde{\Gam}}+d-2,\gamma_{d+1})} \right\rangle_{\J(k+l+\abs{\tilde{\Gam}}+d-2,\gamma_{d+1})}
\end{multline*}
for all $p\in\Pi_{k}^{d-1}$ and $q\in\Pi_{l}^{d-1}$.
\begin{proof}
The result follows directly by applying the change of variables $(x_{1},\dots,x_{d}) = \Psi_{d}(\bd{\zeta},\eta)$ to the integral underlying the left-hand side of the identity.
\end{proof}

\begin{remark}\label{orthoT}
From \autoref{prop:normT} and the orthogonality properties of Jacobi polynomials, we deduce that if $p\in\mathcal{V}_{k}^{\mathrm{ridge},\tilde{\Gam}}$ and $q\in\mathcal{V}_{l}^{\mathrm{ridge},\tilde{\Gam}}$, then $\mathcal{T}^{\Gam}_{n,k}(p)\perp_{\W_{\mathrm{ridge},\Gam}}\mathcal{T}^{\Gam}_{m,l}(q)$ when $(n,k)\neq(m,l)$.
\end{remark}

\begin{remark}\label{rem:injectivity-tensorization}
In view of \autoref{prop:normT} and the fact that the involved Jacobi polynomial is not null (cf.\ \eqref{Jacobi1}), it follows that the operator $\mathcal{T}_{l,k}^{\Gam}$ introduced in \autoref{def:tensorization} is injective.
\end{remark}
\end{proposition}

We close this section with \autoref{descompobreveV} and \autoref{orthoshiftT}, which are orthogonal decomposition results that function as higher-dimensional analogues of \autoref{orthobreve} and \autoref{orthw}, respectively.

\begin{proposition}\label{descompobreveV}
Let $d\in\mathbb{N}_{\geq 3}$, $n\in\mathbb{N}_{0}$, and $\Gam = (\tilde{\Gam},\gamma_{d+1})\in\Lambda^{d}$.
Then,
\begin{equation}\label{decomdepletedspace}
\mathcal{V}_{n}^{\mathrm{ridge},\Gam} = \bigoplus_{k=0}^{n}\mathcal{T}_{n,k}^{\Gam}(\mathcal{V}_{k}^{\mathrm{ridge},\tilde{\Gam}}).
\end{equation}
Moreover, this decomposition is $\W_{\mathrm{ridge},\Gam}$-orthogonal.
\begin{proof}
Let $i,j\in[n]_{0}$ be such that $i\neq j$.
By \autoref{orthoT}, we have that
\begin{equation*}
\mathcal{T}^{\Gam}_{n,i}(\mathcal{V}_{i}^{\mathrm{ridge},\tilde{\Gam}})\perp_{\W_{\mathrm{ridge},\Gam}}\mathcal{T}^{\Gam}_{n,j}(\mathcal{V}_{j}^{\mathrm{ridge},\tilde{\Gam}}).
\end{equation*}
Consequently, the decomposition on the right-hand side of \eqref{decomdepletedspace} is $\W_{\mathrm{ridge},\Gam}$-orthogonal.
This fact, together with \autoref{rem:injectivity-tensorization} and \eqref{dimOPS}, leads to
\begin{multline}\label{dimoplusT}
\operatorname{dim}\mathopen{}\left(\bigoplus_{k=0}^{n}\mathcal{T}^{\Gam}_{n,k}(\mathcal{V}^{\mathrm{ridge},\tilde{\Gam}}_{k})\mathclose{}\right) = \sum_{k=0}^{n}\operatorname{dim}\mathopen{}\left(\mathcal{T}_{n,k}^{\Gam}(\mathcal{V}_{k}^{\mathrm{ridge},\tilde{\Gam}})\mathclose{}\right)\\
= \sum_{k=0}^{n}\binom{k+d-2}{k} = \binom{n+d-1}{n}.
\end{multline}
In view of \eqref{dimOPS} and the above, it follows that
\begin{equation*}
\operatorname{dim}\left(\bigoplus_{k=0}^{n}\mathcal{T}^{\Gam}_{n,k}(\mathcal{V}^{\mathrm{ridge},\tilde{\Gam}}_{k})\right) = \operatorname{dim}\left(\mathcal{V}_{n}^{\mathrm{ridge},\Gam}\right).
\end{equation*}
We shall now show that
\begin{equation}\label{incluoplusT}
\bigoplus_{k=0}^{n}\mathcal{T}^{\Gam}_{n,k}(\mathcal{V}^{\mathrm{ridge},\tilde{\Gam}}_{k})\subset \mathcal{V}_{n}^{\mathrm{ridge},\Gam},
\end{equation}
which, combined with \eqref{dimoplusT}, establishes the equality in \eqref{decomdepletedspace}.
For $m\in\mathbb{N}_{0}$, we first observe that \autoref{orthoT} and \eqref{dimoplusT} imply
\begin{equation*}
\dim\left(\bigoplus_{l=0}^{m}\bigoplus_{r=0}^{l}\mathcal{T}^{\Gam}_{l,r}\left(\mathcal{V}_{r}^{\mathrm{ridge},\tilde{\Gam}}\right)\right) = \sum_{l=0}^{m}\binom{l+d-1}{l} = \binom{m+d}{m} = \operatorname{dim}(\Pi_{m}^{d}).
\end{equation*}
Since $\mathcal{T}^{\Gam}_{l,r}\left(\Pi_{r}^{d-1}\right)\subset\Pi_{l}^{d}$ (cf.~\autoref{def:tensorization}), we deduce that
\begin{equation}\label{descomPoly}
\Pi_{m}^{d} = \bigoplus_{l=0}^{m}\bigoplus_{r=0}^{l}\mathcal{T}^{\Gam}_{l,r}\left(\mathcal{V}_{r}^{\mathrm{ridge},\tilde{\Gam}}\right).
\end{equation}
Thus, given a polynomial $p_{n}\in\bigoplus_{k=0}^{n}\mathcal{T}_{n,k}^{\Gam}( \mathcal{V}_{k}^{\mathrm{ridge},\tilde{\Gam}} )\subset\Pi_{n}^{d}$ and any polynomial $q\in\Pi_{n-1}^{d}$, we can expand the latter according to \eqref{descomPoly} and so we can apply \autoref{orthoT} to conclude that $\langle p_{n},q_{n-1} \rangle_{\mathrm{ridge},\Gam} = 0.$
\end{proof}
\end{proposition}

\begin{proposition}\label{orthoshiftT}
Given $d\in\mathbb{N}_{\geq 3}$ and $n\in\mathbb{N}_{0}$, let $\Gam = (\tilde{\Gam},\gamma_{d+1})\in(-1,\infty)^{d+1}$.
Then, the decomposition
\begin{equation*}
\mathcal{V}_{n}^{\mathrm{ridge},\Gam+e_1+e_2} = \bigoplus_{k=0}^{n}\mathcal{T}_{n,k}^{\Gam+e_{1}+e_{2}}(\mathcal{V}_{k}^{\mathrm{ridge},\tilde{\Gam}+e_{1}+e_{2}})
\end{equation*}
given in \autoref{descompobreveV} is also $\W_{\Gam}$-orthogonal.
\begin{proof}
Let $l,r\in[n]_0$, with $p_{l}\in\mathcal{V}_{l}^{\mathrm{ridge},\tilde{\Gam}+e_{1}+e_{2}}$ and $p_{r}\in\mathcal{V}_{r}^{\mathrm{ridge},\tilde{\Gam}+e_{1}+e_{2}}$.
Applying the change of variables $(x_{1},\dots,x_{d}) = \Psi_{d}(\bd{\zeta},\eta)$ yields the factorization:
\begin{multline}\label{T-inner-products-factorization}
\langle \mathcal{T}_{n,l}^{\Gam+e_{1}+e_{2}}(p_{l}),\mathcal{T}_{n,r}^{\Gam+e_{1}+e_{2}}(p_{r})\rangle_{\Gam} = 2^{-(l+r+\abs{\Gam}+d)}\langle p_{l},p_{r} \rangle_{\tilde{\Gam}}\\
\times\langle \mathrm{P}_{n-l}^{(2l+\abs{\tilde{\Gam}}+d,\gamma_{d+1})},\mathrm{P}_{n-r}^{(2r+\abs{\tilde{\Gam}}+d,\gamma_{d+1})} \rangle_{\mathrm{J}(l+r+\abs{\tilde{\Gam}}+d-1,\gamma_{d+1})}.
\end{multline}
The second term on the right-hand side of the above equality can be written as follows:
\begin{equation}\label{wgortho}
\int_{-1}^{1}\mathrm{P}_{n-l}^{(2l+\abs{\tilde{\Gam}}+d,\gamma_{d+1})}(\eta)\underbrace{\mathrm{P}_{n-r}^{(2r+\abs{\tilde{\Gam}}+d,\gamma_{d+1})}(\eta)(1-\eta)^{r-l-1}}_{Y(\eta)} (1-\eta)^{2l+\abs{\tilde{\Gam}}+d}(1+\eta)^{\gamma_{d+1}}\ \textup{d}\eta.
\end{equation}
If $r\geq l+1$, then $Y\in\Pi_{n-l-1}^{1}$, and since $n\geq r \geq l+1$ the integral appearing in \eqref{wgortho} vanishes, so
\begin{equation*}
\langle \mathcal{T}_{n,l}^{\Gam+e_{1}+e_{2}}(p_{l}),\mathcal{T}_{n,r}^{\Gam+e_{1}+e_{2}}(p_{r})\rangle_{\Gam} = 0.
\end{equation*}
If $l\geq r+1$, the approach is analogous to the above case, by swapping the roles of $r$ and $l$.
\end{proof}
\end{proposition}

\section{Symmetric Group Action on Orthogonal Polynomial Spaces}\label{sec:SGAonOPS}

An important property of the $\LL_{\Gam}^{2}$-orthogonal polynomials on the simplex $\tri^{d}$ is the covariance between the permutations of the parameters $(\gamma_{1},\dots,\gamma_{d+1})$ of the former and the set of vertices $\mathbb{V}_{d} := \{v_{1},\dots,v_{d+1}\}$ of the latter, where
\begin{equation}\label{vertices}
(\forall\,i \in [d]) \quad v_{i} = e_{i} \qquad\text{and}\qquad v_{d+1} = (0,\dotsc,0).
\end{equation}
These properties will be exploited using suitable affine maps which are in one-to-one relation with the permutations of the set $\mathbb{V}_{d}$.

Let $(\mathcal{S}_{d+1},\circ)$ denote the group of permutations on $[d+1]$; i.e., the \emph{symmetric group} of order $d+1$.
Given $\sigma\in\mathcal{S}_{d+1}$, we denote by $T_{\sigma}:\tri^{d}\longrightarrow\tri^{d}$ the affine map
\begin{equation}\label{eq:affine-map-T-sigma}
T_{\sigma}(x) = \sum_{i=1}^{d}x_{i}v_{\sigma(i)} + (1-\abs{x})v_{\sigma(d+1)}
\end{equation}
which satisfies $T_{\sigma}(v_{i}) = v_{\sigma(i)}$ for all $i\in[d+1]$.
Given any function $f$ on $\tri^{d}$, we introduce the pullback notation
\begin{equation}\label{eq:pullback-T-sigma}
T_{\sigma}^{\star}f := f\circ T_{\sigma}.
\end{equation}
Now, we will explore the action of $T_{\sigma}^{\star}$ on the weight function $\W_{\Gam}$ and on the orthogonal polynomial space $\mathcal{V}^{\Gam}_{n}$.
\begin{proposition}\label{prop4}
Let $d\in\mathbb{N}$ and let $\sigma,\tau\in\mathcal{S}_{d+1}$.
Then, $T_{\sigma}\circ T_{\tau} = T_{\sigma\circ\tau}$.
\begin{proof}
Let $x\in\tri^{d}$.
By definition, the affine map $T_{\tau}$ can be written as
\begin{equation*}
T_{\tau}(x) = \sum_{i=1}^{d}x_{i}v_{\tau(i)} + (1-\abs{x})v_{\tau(d+1)},
\end{equation*}
so that $T_{\tau}(x)$ is an affine combination of the vertices $v_{\tau(1)},\dots,v_{\tau(d+1)}$.
Since $T_{\sigma}$ is also an affine map, it preserves affine combinations, and therefore we have
\begin{multline*}
T_{\sigma}\circ T_{\tau}(x) = \sum_{i=1}^{d}x_{i}T_{\sigma}(v_{\tau(i)}) + (1-\abs{x})T_{\sigma}(v_{\tau(d+1)})\\
= \sum_{i=1}^{d}x_{i}v_{\sigma\circ\tau(i)} + (1-\abs{x})v_{\sigma(\tau(d+1))}
= T_{\sigma\circ\tau}(x)
\end{multline*}
as desired.
\end{proof}
\end{proposition}
\begin{proposition}\label{prop5}
Given $d\in\mathbb{N}$ and $\sigma\in\mathcal{S}_{d+1}$, the map $T_{\sigma}$ is a bijection from $\tri^{d}$ onto itself and $T_{\sigma}^{-1} = T_{\sigma^{-1}}$.
\begin{proof}
To show that $T_{\sigma}$ is injective, let $x,y\in\tri^{d}$, such that $T_{\sigma}(x) = T_{\sigma}(y)$. Then,
\begin{equation*}
\sum_{i=1}^{d}(x_{i}-y_{i})v_{\sigma(i)} + \abs{y-x}v_{\sigma(d+1)} = 0
\end{equation*}
if $\sigma(i)\neq d+1$ for all $i\in[d]$, then the set $\{v_{\sigma(i)}\}_{i=1}^{d}$ is linearly independent and $v_{\sigma(d+1)} = 0$, so $x = y$.
On the another hand, if there is one $j\in[d]$ such that $\sigma(j) = d+1$, then the set $\mathbb{V}_{d}\backslash\{v_{\sigma(j)}\}$ is linearly independent.
In this case, we obtain that $x_{i} = y_{i}$ for all $i\in[d]\backslash\{j\}$ and moreover, since $\abs{x} = \abs{y}$, it follows that $x_{j} = y_{j}$.
Thus, in all cases, we conclude that $x = y$.

To show that $T_{\sigma}$ is surjective, let $y\in\tri^{d}$ be given.
As there exists $\sigma^{-1}\in \mathcal{S}_{d+1}$ such that $\sigma\circ \sigma^{-1} = id$, we take $x = T_{\sigma^{-1}}(y)$ and by \autoref{prop4} we have that
\begin{equation*}
T_{\sigma}(x) = T_{\sigma}\circ T_{\sigma^{-1}}(y) = T_{\sigma\circ \sigma^{-1}}(y) = T_{id}(y) = y.
\end{equation*}

To show that $T_{\sigma}^{-1} = T_{\sigma^{-1}}$, it suffices to compute $T_{\sigma}\circ T_{\sigma^{-1}}(x)$ using \autoref{prop4}.
\end{proof}
\end{proposition}
\begin{proposition}\label{prop6}
Given $\Gam\in(-1,\infty)^{d+1}$ and $\sigma\in\mathcal{S}_{d+1}$, the following hold:
\begin{enumerate}[label=(\roman*), ref=(\roman*)]
\setcounter{enumi}{0}
\item\label{prop6one} $T_{\sigma}^{\star}\W_{\Gam} = \W_{\sigma(\Gam)}$, where $\sigma(\Gam) := (\gamma_{\sigma(1)},\dots,\gamma_{\sigma(d+1)})$.
\item\label{prop6onehalf} For all $f,g\in \LL^2_{\Gam}$, we have that $\langle T_{\sigma}^{\star}f,T_{\sigma}^{\star}g\rangle_{\sigma(\Gam)} = \langle f,g\rangle_{\Gam}$.
\item\label{prop6two} $T_{\sigma}^{\star}\mathcal{V}_{n}^{\Gam} = \mathcal{V}_{n}^{\sigma(\Gam)}$.
\end{enumerate}
\begin{proof}
To prove part \ref{prop6one}, we introduce the coordinates $\xi(x) = (x_{1},\dots,x_{d},1-\abs{x})$; i.e., the \emph{barycentric coordinates} of $x$ with respect to the vertices of $\tri^{d}$.
Then,
\begin{equation*}
T_{\sigma}(x) = \sum_{i=1}^{d+1}\xi_{i}(x)v_{\sigma(i)} \overset{j = \sigma(i)}{=} \sum_{j=1}^{d+1}\xi_{\sigma^{-1}(j)}(x)v_{j}.
\end{equation*}
It follows that $\xi(T_{\sigma}(x)) = (\xi_{\sigma^{-1}(1)}(x),\dots,\xi_{\sigma^{-1}(d+1)}(x))$.
Using this fact and the fact that $\W_{\Gam}(x) \stackrel{\eqref{weight}}{=} \prod_{i=1}^{d+1} \xi_{i}(x)^{\gamma_i}$, we find that
\begin{equation*}
T_{\sigma}^{\star}\W_{\Gam}(x) = \prod_{i=1}^{d+1} \xi_{i}(T_{\sigma}(x))^{\gamma_{i}} = \prod_{i=1}^{d+1} \xi_{\sigma^{-1}(i)}(x)^{\gamma_{i}} \overset{i = \sigma(j)}{=} \prod_{j=1}^{d+1} \xi_{j}(x)^{\gamma_{\sigma(j)}} = \W_{\sigma(\Gam)}(x).
\end{equation*}
On the other hand, we observe that
\begin{equation*}
\langle T_{\sigma}^{\star}f,T_{\sigma}^{\star}g\rangle_{\sigma(\Gam)} = \int_{\tri^{d}}f(T_{\sigma}(x))g(T_{\sigma}(x))\W_{\sigma(\Gam)}(x)\ \textup{d}x.
\end{equation*}
Performing the change of variable $z = T_{\sigma}(x)$, and using the fact that $\abs{\operatorname{det}(\operatorname{D}T_{\sigma})} = 1$ and \autoref{prop5}, we obtain
\begin{equation*}
\langle T_{\sigma}^{\star}f,T_{\sigma}^{\star}g\rangle_{\sigma(\Gam)} = \int_{\tri^{d}}f(z)g(z)T_{\sigma^{-1}}^{\star}\W_{\sigma(\Gam)}(z)\ \textup{d}z \stackrel{\text{\ref{prop6one}}}{=} \langle f,g\rangle_{\Gam}.
\end{equation*}
For part \ref{prop6two}, let $p\in\mathcal{V}_{n}^{\Gam}$.
Given any $q\in\Pi_{n-1}^{d}$, using part \ref{prop6onehalf} and the fact that $T_{\sigma^{-1}}^{\star}$ preserves the total degree of polynomials,
\begin{equation*}
\langle T_{\sigma}^{\star}p,q \rangle_{\sigma(\Gam)} = \langle p, T^{\star}_{\sigma^{-1}}q\rangle_{\Gam} = 0.
\end{equation*}
This shows that $T_{\sigma}^{\star}\mathcal{V}_{n}^{\Gam}\subset\mathcal{V}_{n}^{\sigma(\Gam)}$.
As $T_{\sigma}^{\star}$ is and injective linear map between $\mathcal{V}_{n}^{\Gam}$ and $\mathcal{V}_{n}^{\sigma(\Gam)}$, because $T_{\sigma}$ is a bijection, we have that $\dim\mathcal{V}_{n}^{\Gam} = \dim\mathcal{V}_{n}^{\sigma(\Gam)}$, so $T_{\sigma}^{\star}\mathcal{V}_{n}^{\Gam} = \mathcal{V}_{n}^{\sigma(\Gam)}$.
\end{proof}
\end{proposition}
\begin{proposition}\label{prop7}
Let $\sigma\in\mathcal{S}_{d+1}$ and $k\in[d]$.
Then:
\begin{enumerate}[label=(\roman*), ref=(\roman*)]
\setcounter{enumi}{0}
\item\label{prop7one} $\partial_{k}T_{\sigma} = v_{\sigma(k)} - v_{\sigma(d+1)}$.
\item\label{prop7two} For all $f\in\CC^{1}(\overline{\tri^{d}})$, we have that
\begin{equation*}
\partial_{k}T_{\sigma}^{\star}f = \sum_{l=1}^{d} (v_{\sigma(k)}-v_{\sigma(d+1)})_{l} \, T_{\sigma}^{\star}\partial_{l}f.
\end{equation*}
\end{enumerate}
\begin{proof}
For part \ref{prop7one}, we have that
\begin{multline*}
\partial_{k}T_{\sigma}(x) = \sum_{i=1}^{d}\partial_{k}x_{i}v_{\sigma(i)} + \partial_{k}(1-\abs{x})v_{\sigma(d+1)}\\
= \sum_{i=1}^{d}\delta_{ik}v_{\sigma(i)} - v_{\sigma(d+1)} = v_{\sigma(k)} - v_{\sigma(d+1)}.
\end{multline*}
On the other hand, for part \ref{prop7two}, using the chain rule and \ref{prop7one} we have
\begin{equation*}
\partial_{k}T_{\sigma}^{\star}f(x) = \sum_{l=1}^{d}T_{\sigma}^{\star}\partial_{l}f(x) \, \partial_{k}(T_{\sigma})_{l}(x) = \sum_{l=1}^{d} (v_{\sigma(k)}-v_{\sigma(d+1)})_{l} \, T_{\sigma}^{\star}\partial_{l}f(x),
\end{equation*}
which completes the proof.
\end{proof}
\end{proposition}
In \autoref{sec:markovinequality}, we will find it convenient to have, for each $i \in [d]$, an affine map $\sigma_{i}\in\mathcal{S}_{d+1}$ such that either $\partial_{1,2}T_{\sigma_{i}}^{\star}f = T_{\sigma_{i}}^{\star}\partial_{i}f$ or $\partial_{1,2}T_{\sigma_{i}}^{\star}f = -T_{\sigma_{i}}^{\star}\partial_{i}f$.
To this end, we construct the subset $\Sigma_{d} \subset \mathcal{S}_{d+1}$ as follows.
\begin{definition}\label{def:Sigmad}
Let $d\in\mathbb{N}$ with $d \geq 2$.
We define $\Sigma_d := \{ \sigma_{1}, \sigma_{2}, \dotsc, \sigma_{d} \}$, where
\begin{enumerate}
\item $\sigma_{1}\in\mathcal{S}_{d+1}$ is defined by $\sigma_{1}(1) = 2$, $\sigma_{1}(2) = d+1$, $\sigma_{1}(d+1) = 1$ and $\sigma_{1}(k) = k$ for all $k\in\{3,\dots,d\}$;
\item $\sigma_{2}\in\mathcal{S}_{d+1}$ is defined by $\sigma_{2}(1) = d+1$, $\sigma_{2}(2) = 1$, $\sigma_{2}(d+1) = 2$ and $\sigma_{2}(k) = k$ for all $k\in\{3,\dots,d\}$;
\item for $i\in\{3,\dots,d\}$, $\sigma_{i}\in\mathcal{S}_{d+1}$ is defined by $\sigma_{i}(1) = d+1$, $\sigma_{i}(2) = i$, $\sigma_{i}(i) = 2$, $\sigma_{i}(d+1) = 1$ and $\sigma_{i}(k) = k$ for all $k\in\{3,\dots,d\}\backslash\{i\}$.
\end{enumerate}
\end{definition}
\begin{proposition}\label{prop8}
Let $f\in\CC^{1}(\overline{\tri^{d}})$ and $i\in\{3,\dots,d\}$.
Then, the following hold:
\begin{enumerate}[label=(\roman*), ref=(\roman*)]
\setcounter{enumi}{0}
\item\label{prop8one} $\partial_{1}T_{\sigma_{1}}^{\star}f = T_{\sigma_{1}}^{\star}\partial_{1,2}f$,
$\partial_{2}T_{\sigma_{1}}^{\star}f = -T_{\sigma_{1}}^{\star}\partial_{1}f$ and $\partial_{1,2}T_{\sigma_{1}}^{\star}f = -T_{\sigma_{1}}^{\star}\partial_{2}f.$
\item\label{prop8two} $\partial_{1}T_{\sigma_{2}}^{\star}f = -T_{\sigma_{2}}^{\star}\partial_{2}f, \quad
\partial_{2}T_{\sigma_{2}}^{\star}f = -T_{\sigma_{2}}^{\star}\partial_{1,2}f$ and $\partial_{1,2}T_{\sigma_{2}}^{\star}f = T_{\sigma_{2}}^{\star}\partial_{1}f.$
\item\label{prop8three} $\partial_{1}T_{\sigma_{i}}^{\star}f = -T_{\sigma_{i}}^{\star}\partial_{1}f, \quad
\partial_{2}T_{\sigma_{i}}^{\star}f = T_{\sigma_{i}}^{\star}\partial_{1,i}f$ and $\partial_{1,2}T_{\sigma_{i}}^{\star}f = T_{\sigma_{i}}^{\star}\partial_{i}f.$
\end{enumerate}
\begin{proof}
Direct from \ref{prop7two} in \autoref{prop7} and \autoref{def:Sigmad}.
\end{proof}
\end{proposition}

\section{Schur-type Inequalities on the Simplex}\label{sec:Schur-type-inequalities}

\subsection{Schur-type inequalities with a single parameter shift}\label{subsec:natural-Schur}
In this subsection we prove the Schur-type inequality \eqref{crude-natural-Schur}.
First, we will proceed by proving a stronger bound which is valid for $\W_{\Gam+e_i}$-orthogonal polynomials only;
this special case is obtained under the additional assumption $i = d+1$ using a suitable integration by parts inspired by \cite[Prop.~3.4(i)]{Figueroa:2017b} and extended to other $i$ by using the covariance result in \autoref{prop6}.
Finally, the desired general result, which is valid for any polynomial, is built upon the previous one with the help of the orthogonal decomposition \eqref{eq:orthogonal-decomposition-of-polynomials} and is encapsulated in \autoref{thm:whole-natural-Schur}.

\begin{lemma}\label{lem:slice-natural-Schur}
Let $d\in\mathbb{N}_{\geq 2}$ and $\Gam\in(-1,\infty)^{d+1}$.
Then, for all $p_{n}\in\mathcal{V}^{\Gam+e_{i}}_{n}$
\begin{equation*}
\norm{p_{n}}_{\Gam}^{2} = \frac{(2n+\abs{\Gam}+d+1)}{\gamma_{i}+1}\norm{p_{n}}_{\Gam+e_{i}}^{2}.
\end{equation*}
\begin{proof}
Let $r\in\mathcal{V}_{n}^{\Gam+e_{d+1}}$.
Then, there exists a homogeneous polynomial $h_{r}$ of degree $n$ such that $r-h_{r}\in\Pi_{n-1}^{d}$ and hence $x\cdot\nabla r-x\cdot\nabla h_{r}\in\Pi_{n-1}^{d}$.
Thus,
\begin{equation}\label{normident1}
\langle x\cdot \nabla r,r \rangle_{\Gam+e_{d+1}} =\langle x\cdot \nabla h_{r},r \rangle_{\Gam+e_{d+1}} = n\langle h_{r},r\rangle_{\Gam+e_{d+1}} = n\norm{r}^{2}_{\Gam+e_{d+1}}.
\end{equation}
Using the fact that $\operatorname{div}(x) = d$, we have
\begin{multline}\label{normident2}
(2n+d)\norm{r}_{\Gam+e_{d+1}}^{2} = 2\langle x\cdot \nabla r,r \rangle_{\Gam+e_{d+1}} + \operatorname{div}(x)\norm{r}_{\Gam+e_{d+1}}^{2}\\
= \int_{\tri^{d}}\operatorname{div}\left(r^{2}(x)x\right)\W_{\Gam+e_{d+1}}(x)\ \textup{d}x.
\end{multline}
Integrating by parts and using that
\begin{equation*}
\left(\forall\,x\in\tri^{d}\right) \quad x\cdot\nabla \W_{\Gam+e_{d+1}}(x) = (\abs{\Gam}+1)\W_{\Gam+e_{d+1}}(x) - (\gamma_{d+1}+1)\W_{\Gam}(x),
\end{equation*}
we can recast \eqref{normident2} as
\begin{equation*}
(2n+d)\norm{r}_{\Gam+e_{d+1}}^{2} = (\gamma_{d+1}+1)\norm{r}^{2}_{\Gam} - (\abs{\Gam}+1)\norm{r}_{\Gam+e_{d+1}}^{2},
\end{equation*}
so
\begin{equation}\label{particularcase}
\norm{r}_{\Gam}^{2} = \frac{(2n+\abs{\Gam}+d+1)}{\gamma_{d+1}+1}\norm{r}_{\Gam+e_{d+1}}^{2},
\end{equation}
which is the desired result in the special case $i=d+1$.

Now, given $i\in[d]$, let $\tau_{i}$ be the permutation in $\mathcal{S}_{d+1}$ satisfying $\tau_{i}(d+1) = i$, $\tau_{i}(i) = d+1$ and $\tau_{i}(k)=k$ for all $k\in[d+1]\backslash\{i,d+1\}$.
The associated affine map (cf.\ \eqref{eq:affine-map-T-sigma}) is
\begin{equation*}
\left(\forall x\in\tri^{d}\right)\ T_{\tau_{i}}(x) = \sum_{k=1}^{d}x_{k}v_{\tau_{i}(k)} + (1-\abs{x})v_{\tau_{i}(d+1)}.
\end{equation*}
Let $p\in\mathcal{V}_{n}^{\Gam + e_{i}}$.
From \ref{prop6two} of \autoref{prop6}, we have that $T_{\tau_{i}}^{\star}p \in \mathcal{V}_{n}^{\tau_{i}(\Gam+e_{i})}$.
As $\mathcal{V}_{n}^{\tau_{i}(\Gam+e_{i})} = \mathcal{V}_{n}^{\tau_{i}(\Gam)+e_{d+1}}$, we can then apply the identity \eqref{particularcase} to obtain
\begin{equation*}
\norm{T_{\tau_{i}}^{\star}p}_{\tau_{i}(\Gam)}^{2} = \frac{(2n+\abs{\tau_{i}(\Gam)}+d+1)}{(\tau_{i}(\Gam))_{d+1}+1}\norm{T_{\tau_{i}}^{\star}p}_{\tau_{i}(\Gam+e_{i})}^{2}.
\end{equation*}
Using part~\ref{prop6onehalf} of \autoref{prop6} to transform the norms and the fact that $\abs{\tau_{i}(\Gam)} = \abs{\Gam}$ and $(\tau_{i}(\Gam))_{d+1} = \gamma_i$, this identity simplifies to the desired result.
\end{proof}
\end{lemma}

\begin{thm}\label{thm:whole-natural-Schur}
Let $d\in\mathbb{N}_{\geq 2}$ and $\Gam\in(-1,\infty)^{d+1}$.
Let $n\in\mathbb{N}_{0}$ and $i\in[d+1]$.
Then, for all $q\in\Pi_{n}^{d}$
\begin{equation*}
\norm{q}_{\Gam}^{2}\leq \frac{(n+1)(n+\abs{\Gam}+d+1)}{\gamma_{i}+1} \norm{q}_{\Gam+e_{i}}^{2}.
\end{equation*}
\begin{proof}
Let $q\in\Pi_{n}^{d}$.
By \eqref{eq:orthogonal-decomposition-of-polynomials}, $\Pi_{n}^{d}=\bigoplus_{k=0}^{n}\mathcal{V}_{k}^{\Gam+e_{i}}$, so via the triangle inequality, we have the estimate
\begin{equation*}
\norm{q}_{\Gam}^{2} \leq \left[ \sum_{k=0}^{n}\norm{\proj_{k}^{\Gam+e_{i}}(q)}_{\Gam} \right]^{2}.
\end{equation*}
On the other hand, from \autoref{lem:slice-natural-Schur} we get the estimate
\begin{equation*}
\norm{\proj_{k}^{\Gam+e_{i}}(q)}_{\Gam} \leq \sqrt{\frac{2k+\abs{\Gam}+d+1}{\gamma_{i}+1}}\norm{\proj_{k}^{\Gam+e_{i}}(q)}_{\Gam+e_{i}}.
\end{equation*}
Thus combining the above two estimates with the Cauchy--Schwarz inequality, we obtain
\begin{multline*}
\norm{q}_{\Gam}^{2} \leq \left( \sum_{k=0}^{n}\frac{2k+\abs{\Gam}+d+1}{\gamma_{i}+1} \right)\left(\sum_{k=0}^{n}\norm{\proj_{k}^{\Gam+e_{i}}(q)}_{\Gam+e_{i}}^{2}\right)\\
= \frac{(n+1)\left( n+ \abs{\Gam}+d+1\right)}{\gamma_{i}+1}\sum_{k=0}^{n}\norm{\proj_{k}^{\Gam+e_{i}}(q)}_{\Gam+e_{i}}^{2}.
\end{multline*}
From here the desired result follows by Parseval's identity.
\end{proof}
\end{thm}

\subsection{Schur-type inequalities with two parameter shifts}
In this subsection we prove the Schur-type inequality \eqref{crude-middle-Schur}.
Again we start with a stronger bound valid for orthogonal polynomials only and here we will obtain it by exploing the results of \autoref{subsec:natural-Schur} in \autoref{lem:slice-middle-Schur}; we will also show that the obtained growth rate of this bound with respect to the polynomial degree is sharp in \autoref{pro:slice-middle-Schur-is-sharp}.
Finally, the desired general result is built upon \autoref{lem:slice-middle-Schur} and \eqref{eq:orthogonal-decomposition-of-polynomials} and is encapsulated in \autoref{thm:whole-middle-Schur}.

\begin{lemma}\label{lem:slice-middle-Schur}
Let $d\in\mathbb{N}_{\geq 2}$, $\Gam\in(-1,\infty)^{d+1}$, $n\in\mathbb{N}_{0}$ and $i,j\in[d+1]$ with $i\neq j$.
Then, for all $p\in\mathcal{V}_{n}^{\Gam+e_{i}+e_{j}}$ we have that
\begin{equation*}
\norm{p}_{\Gam} \leq\sqrt{\frac{(n+1)(n+\abs{\Gam}+d+1)(2n+\abs{\Gam}+d+2)}{(\gamma_{i}+1)(\gamma_{j}+1)}} \norm{p}_{\Gam+e_{i}+e_{j}}.
\end{equation*}
\begin{proof}
Let $p\in\mathcal{V}_{n}^{\Gam+e_{i}+e_{j}}$.
Then,
\begin{multline*}
\norm{p}_{\Gam} \overset{\text{Th.\,\ref{thm:whole-natural-Schur}}}{\leq} \left(\frac{(n+1)\left( n+ \abs{\Gam}+d+1\right)}{\gamma_{i}+1}\right)\norm{p}_{\Gam+e_{i}}^{2}\\
\overset{\text{Lem.\,\ref{lem:slice-natural-Schur}}}{=} \left(\frac{(n+1)\left( n+ \abs{\Gam}+d+1\right)}{\gamma_{i}+1}\right)\frac{(2n+\abs{\Gam}+d+2)}{\gamma_{j}+1}\norm{p}_{\Gam+e_{i}+e_{j}}^{2}.
\end{multline*}
\end{proof}
\end{lemma}

\begin{proposition}\label{pro:slice-middle-Schur-is-sharp}
Let $d\in\mathbb{N}_{\geq 2}$, $\Gam\in(-1,\infty)^{d+1}$ and $i,j\in[d+1]$ with $i\neq j$.
Suppose that, for some $C > 0$ and some $n_0 \in \mathbb{N}$ there holds
\begin{equation*}
(\forall\,n\geq n_0)\ (\forall\,p \in\mathcal{V}_{n}^{\Gam+e_i+e_j}) \quad
\norm{p}_{\Gam} \leq C n^r \norm{p}_{\Gam+e_i+e_j}.
\end{equation*}
Then, $r \geq 3/2$.
In this sense, \autoref{lem:slice-middle-Schur} is sharp.
\begin{proof}
Let us recall the notation $\Gam := (\tilde{\Gam},\gamma_{d+1})$.
Given $n \in \mathbb{N}$, let $Q_{n}^{\Gam+e_{1}+e_{2}}$ be the polynomial defined by
\begin{equation}\label{eq:Bad-guy}
Q_{n}^{\Gam+e_{1}+e_{2}}(x) = \PP_{n}^{(\abs{\tilde{\Gam}}+d+1,\gamma_{d+1})}(1-2\abs{x}).
\end{equation}
If $d = 2$, we can write \eqref{eq:Bad-guy} as
\begin{equation*}
Q_{n}^{\Gam+e_{1}+e_{2}} = T_{\sigma_{1}}^{\star}P_{0,n}^{\gamma_3,\gamma_1+1,\gamma_2+1},
\end{equation*}
where $P_{0,n}^{\gamma_3,\gamma_1+1,\gamma_2+1}$ is defined and proved to belong to $\mathcal{V}_{n}^{(\gamma_3,\gamma_1+1,\gamma_2+1)}$ in \cite[Prop.~2.4.1]{DunklXu:2014}.
Then, by part~\ref{prop6two} of \autoref{prop6}, we have $Q_{n}^{\Gam+e_{1}+e_{2}}\in\mathcal{V}_{n}^{\Gam+e_{1}+e_{2}}$.
Next, performing the change of variables $x = \Psi_{2}(\zeta,\eta)$ (cf.\ \eqref{Psi}) on the integrals underlying the $\W_{\Gam}$- and $\W_{\Gam+e_1+e_2}$-weighted norms and applying \autoref{prop:double-shift}, we obtain
\begin{equation*}
\frac{\norm{Q_{n}^{\Gam+e_{1}+e_{2}}}_{\Gam}^{2}}{\norm{Q_{n}^{\Gam+e_{1}+e_{2}}}_{\Gam+e_{1}+e_{2}}^{2}}
= 16 \, \frac{\norm{1}_{\J(\gamma_{1},\gamma_{2})}^{2}}{\norm{1}_{\J(\gamma_{1}+1,\gamma_{2}+1)}^{2}} \, \mathfrak{C}^{(\gamma_{1}+\gamma_{2}+1,\gamma_{3})}(n).
\end{equation*}
If $d \geq 3$, it follows from \autoref{def:tensor-mathbb-T} that $Q_{n}^{\Gam+e_{1}+e_{2}} = \mathbb{T}^{\Gam+e_{1}+e_{2}}_{n,0}(1)$.
From \autoref{prop:decomposition-V}, it follows that $Q_{n}^{\Gam+e_{1}+e_{2}}\in\mathcal{V}_{n}^{\Gam+e_{1}+e_{2}}$.
Then, applying \autoref{prop:poly-tensor} and \autoref{prop:double-shift}, we obtain
\begin{equation*}
\frac{\norm{Q_{n}^{\Gam+e_{1}+e_{2}}}_{\Gam}^{2}}{\norm{Q_{n}^{\Gam+e_{1}+e_{2}}}_{\Gam+e_{1}+e_{2}}^{2}}
= 4 \, \frac{\norm{1}_{\tilde{\Gam}}^{2}}{\norm{1}_{\tilde{\Gam}+e_{1}+e_{2}}^{2}} \, \mathfrak{C}^{(\abs{\Tilde{\Gam}}+d-1,\gamma_{d+1})}(n).
\end{equation*}
In both cases, since $\mathfrak{C}^{(\abs{\Tilde{\Gam}}+d-1,\gamma_{d+1})}$ is a polynomial of degree three in $n$, the limit
\begin{equation*}
\lim_{n\to\infty}\frac{\norm{Q_{n}^{\Gam+e_{1}+e_{2}}}_{\Gam}^{2}}{n^{3}\norm{Q_{n}^{\Gam+e_{1}+e_{2}}}_{\Gam+e_{1}+e_{2}}^{2}}
\end{equation*}
exists and equals a positive constant.

Next, let $\tau\in\mathcal{S}_{d+1}$ be any permutation satisfying $\tau(i) = 1$ and $\tau(j) = 2$.
From \ref{prop6two} of \autoref{prop6}, we find that $T_{\tau}^{\star}Q_{n}^{\tau^{-1}(\Gam)+e_{1}+e_{2}} = T_{\tau}^{\star}Q_{n}^{\tau^{-1}(\Gam+e_{i}+e_{j})} \in \mathcal{V}_{n}^{\Gam+e_{i}+e_{j}}$.
By virtue of \ref{prop6onehalf} of \autoref{prop6}, we have that
\begin{equation*}
\frac{\norm{T_{\tau}^{\star}Q_{n}^{\tau^{-1}(\Gam)+e_{1}+e_{2}}}_{\tau\circ\tau^{-1}(\Gam)}^{2}}{\norm{T_{\tau}^{\star}Q_{n}^{\tau^{-1}(\Gam)+e_{1}+e_{2}}}_{\tau(\tau^{-1}(\Gam)+e_{1}+e_{2})}^{2}}
= \frac{\norm{Q_{n}^{\tau^{-1}(\Gam)+e_{1}+e_{2}}}_{\tau^{-1}(\Gam)}^{2}}{\norm{Q_{n}^{\tau^{-1}(\Gam)+e_{1}+e_{2}}}_{\tau^{-1}(\Gam)+e_{1}+e_{2}}^{2}}
\end{equation*}
This relation, combined with the previous limit, implies that there exists a constant $\tilde{C}>0$ such that
\begin{equation*}
\lim_{n\to \infty}\frac{\norm{T_{\tau}^{\star}Q_{n}^{\tau^{-1}(\Gam)+e_{1}+e_{2}}}_{\Gam}}{n^{3/2}\norm{T_{\tau}^{\star}Q_{n}^{\tau^{-1}(\Gam)+e_{1}+e_{2}}}_{\Gam+e_{i}+e_{j}}}
= \tilde{C}.
\end{equation*}
Suppose that for some power $r < 3/2$ there existed $n_0 \in \mathbb{N}_0$ and $C>0$ such that $\norm{p}_{\Gam} \leq C n^{r} \norm{p}_{\Gam+e_i+e_j}$ for all $p \in \mathcal{V}_{n}^{\Gam+e_i+e_j}$.
Then, the above limit would vanish, which contradicts the fact that $\tilde C$ is positive.
\end{proof}
\end{proposition}

\begin{thm}\label{thm:whole-middle-Schur}
Let $d\in\mathbb{N}_{\geq2}$, $\Gam\in(-1,\infty)^{d+1}$, $n\in\mathbb{N}_{0}$ and $i,j\in[d+1]$ with $i \neq j$.
Then, for all $p\in\Pi_{n}^{d}$ we have that
\begin{equation*}
\norm{p}_{\Gam}\leq \sqrt{\frac{(n+1)(n+2)(n+\abs{\Gam}+d+1)(n+\abs{\Gam}+d+2)}{2(\gamma_{i}+1)(\gamma_{j}+1)}}\norm{p}_{\Gam+e_{i}+e_{j}}.
\end{equation*}
\begin{proof}
Let $p\in\Pi_{n}^{d}$.
From the decomposition $\Pi_{n}^{d} = \bigoplus_{k=0}^{n}\mathcal{V}_{k}^{\Gam+e_{i}+e_{j}}$, the triangle inequality and \autoref{lem:slice-middle-Schur},
\begin{equation*}
\norm{p}_{\Gam}^{2} \leq \left(\sum_{k=0}^{n}\norm{\proj_{k}^{\Gam+e_{i}+e_{j}}(p)}_{\Gam}\right)^{2}
\leq \left(\sum_{k=0}^{n} \sqrt{\mathcal{D}^{\Gam}(k)} \norm{\proj_{k}^{\Gam+e_{i}+e_{j}}(p)}_{\Gam+e_i+e_j}\right)^{2},
\end{equation*}
where
\begin{equation*}
\mathfrak{D}^{\Gam}(k) := \frac{(k+1)(k+\abs{\Gam}+d+1)(2k+\abs{\Gam}+d+2)}{(\gamma_{i}+1)(\gamma_{j}+1)}.
\end{equation*}
Hence, via the Cauchy--Schwarz inequality in $\mathbb{R}^{n+1}$ and Parseval's identity we obtain
\begin{equation*}
\norm{p}_{\Gam}^{2} \leq \left(\sum_{k=0}^{n}\mathfrak{D}^{\Gam}(k)\right)\norm{p}_{\Gam+e_{i}+e_{j}}^{2}.
\end{equation*}
The desired result then follows from the fact that
\begin{equation*}
\sum_{k=0}^{n}\mathfrak{D}^{\Gam}(k) = \frac{(n+1)(n+2)(n+\abs{\Gam}+d+1)(n+\abs{\Gam}+d+2)}{2(\gamma_{i}+1)(\gamma_{j}+1)}.
\end{equation*}
\end{proof}
\end{thm}

\subsection{Schur-type inequalities with two parameter shifts and the addition of a ridge singularity}
In this subsection we prove the Schur-type inequality \eqref{crude-depleted-Schur}.
Once again we start with a stronger bound valid for orthogonal polynomials only, this time through an induction argument on the dimension of the simplex exploiting results of \autoref{sec:Jacpoly} and \autoref{subsec:Dep-OP}, in \autoref{lem:slice-depleted-Schur}.
Then, the desired general result is built upon \autoref{lem:slice-depleted-Schur} and \eqref{eq:orthogonal-decomposition-of-polynomials} and is encapsulated in \autoref{thm:whole-depleted-Schur}.

\begin{lemma}\label{lem:slice-depleted-Schur}
Let $d\in\mathbb{N}_{\geq 2}$, $\Gam\in(-1,\infty)^{d+1}$ and $n\in\mathbb{N}_{0}$.
Then,
\begin{equation*}
\left(q\in\mathcal{V}_{n}^{\mathrm{ridge},\Gam+e_{1}+e_{2}}\right)~\norm{q}_{\Gam}^{2} \leq \frac{(\gamma_{1}+\gamma_{2}+3)(2n+\abs{\Gam}+d+1)}{(\gamma_{1}+1)(\gamma_{2}+1)}\norm{q}_{\mathrm{ridge},\Gam+e_{1}+e_{2}}^{2}.
\end{equation*}
\begin{proof}
We proceed by induction on the dimension.

\textbf{The base case}.
Let $n\in\mathbb{N}_{0}$ and $j\in[n]_{0}$.
Using \eqref{eq:decom-depled-weight-2d} from the proof of \autoref{orthobreve} (with $\gamma_1$ and $\gamma_2$ shifted) and \eqref{eq:decom-inner-product-2d} from the proof of \autoref{orthw}, we find that
\begin{multline*}
\frac{\norm{\PP_{n,j}^{\mathrm{ridge},\Gam+e_{1}+e_{2}}}_{\Gam}^{2}}{\norm{\PP_{n,j}^{\mathrm{ridge},\Gam+e_{1}+e_{2}}}_{\mathrm{ridge},\Gam+e_{1}+e_{2}}^{2}}
= 2^{3} \,
\frac{\norm{\PP_{j}^{(\gamma_{1}+1,\gamma_{2}+1)}}_{\J(\gamma_{1},\gamma_{2})}^{2}}{\norm{\PP_{j}^{(\gamma_{1}+1,\gamma_{2}+1)}}_{\J(\gamma_{1}+1,\gamma_{2}+1)}^{2}}\\
\times \frac{\norm{\PP_{n-j}^{(2j+\gamma_{1}+\gamma_{2}+2,\gamma_{3})}}_{\J(2j+\gamma_{1}+\gamma_{2}+1,\gamma_{3})}^{2}}{\norm{\PP_{n-j}^{(2j+\gamma_{1}+\gamma_{2}+2,\gamma_{3})}}_{\J(2j+\gamma_{1}+\gamma_{2}+2,\gamma_{3})}^{2}}.
\end{multline*}
By evaluating the norm ratios of Jacobi polynomials above with the help of \autoref{parshift} and \autoref{totalshift}, we obtain
\begin{multline*}
\norm{\PP_{n,j}^{\mathrm{ridge},\Gam+e_{1}+e_{2}}}_{\Gam}^{2}
= \frac{(2j+\gamma_{1}+\gamma_{2}+3)(\gamma_{1}+\gamma_{2}+2)}{(\gamma_{1}+1)(\gamma_{2}+1)}\\
\times \frac{(2n+\gamma_{1}+\gamma_{2}+\gamma_{3}+3)}{(2j+\gamma_{1}+\gamma_{2}+2)} \norm{\PP_{n,j}^{\mathrm{ridge},\Gam+e_{1}+e_{2}}}_{\mathrm{ridge},\Gam+e_{1}+e_{2}}^{2}.
\end{multline*}
The rational expression involving $n$, $j$ and the entries of $\Gam$ on the right-hand side of the above equation attains its maximum among the $j$ in $[n]_{0}$ at $j=0$, so
\begin{equation*}
\norm{\PP_{n,j}^{\mathrm{ridge},\Gam+e_{1}+e_{2}}}_{\Gam}^{2} \leq \frac{(\gamma_{1}+\gamma_{2}+3)(2n+\abs{\Gam}+3)}{(\gamma_{1}+1)(\gamma_{2}+1)}\norm{\PP_{n,j}^{\mathrm{ridge},\Gam+e_{1}+e_{2}}}_{\mathrm{ridge},\Gam+e_{1}+e_{2}}^{2}.
\end{equation*}
As by \autoref{orthobreve} and \autoref{orthw} the set $\{ \PP_{n,j}^{\mathrm{ridge},\Gam+e_{1}+e_{2}} \}_{j=0}^n$ is simultaneously a $\W_{\mathrm{ridge},\Gam+e_1+e_2}$- and $\W_{\Gam}$-orthogonal basis of $\mathcal{V}_{n}^{\mathrm{ridge},\Gam+e_1+e_2}$, the above inequality generalizes to
\begin{equation*}
(\forall \, p_{n} \in \mathcal{V}_{n}^{\mathrm{ridge},\Gam+e_{1}+e_{2}}) \quad
\norm{p_{n}}^{2}_{\Gam}\leq \frac{(\gamma_{1}+\gamma_{2}+3)(2n+\abs{\Gam}+3)}{(\gamma_{1}+1)(\gamma_{2}+1)}\norm{p_{n}}^{2}_{\mathrm{ridge},\Gam+e_{1}+e_{2}},
\end{equation*}
which is the desired result in the case $d = 2$.

\textbf{The induction step.}
Let us assume that the desired result is valid for a dimension $d \geq 2$ and let us prove it for the dimension $d+1$.
Thus, let $\Gam = (\tilde{\Gam},\gamma_{d+2})\in(-1,\infty)^{d+2}$, $n\in\mathbb{N}_0$, and $p\in\mathcal{V}_{n}^{\mathrm{ridge},\Gam+e_{1}+e_{2}}$.
By \autoref{descompobreveV} there exist $q_{k}\in\mathcal{V}_{k}^{\mathrm{ridge},\tilde{\Gam}+e_{1}+e_{2}}$, $k \in [n]_0$, such that
\begin{equation*}
p = \sum_{k=0}^{n} \mathcal{T}_{n,k}^{\Gam+e_{1}+e_{2}}(q_{k}).
\end{equation*}
In view of \autoref{orthoshiftT}, we have
\begin{equation}\label{eq:normrn}
\norm{p}_{\Gam}^{2} = \sum_{k=0}^{n} \norm{\mathcal{T}_{n,k}^{\Gam+e_{1}+e_{2}}(q_{k})}_{\Gam}^{2}.
\end{equation}
Given $k \in [n]_0$, using \eqref{T-inner-products-factorization} from the proof of \autoref{orthoshiftT} yields
\begin{equation*}
\norm{\mathcal{T}_{n,k}^{\Gam+e_{1}+e_{2}}(q_{k})}_{\Gam}^{2} = 2^{-(2k+\abs{\Gam}+d+1)}\norm{q_{k}}_{\tilde{\Gam}}^{2}
\norm{\PP_{n-k}^{\left(2k+\abs{\tilde{\Gam}}+d+1,\gamma_{d+2}\right)}}_{\J\left(2k+\abs{\tilde{\Gam}}+d,\gamma_{d+2}\right)}^{2}.
\end{equation*}
By applying \autoref{parshift} on the last factor on the right-hand side of above equality, we obtain
\begin{multline}\label{eq:identnormT}
\norm{\mathcal{T}_{n,k}^{\Gam+e_{1}+e_{2}}(q_{k})}_{\Gam}^{2} = 2^{-(2k+\abs{\Gam}+d+1)}\norm{q_{k}}_{\tilde{\Gam}}^{2}\\
\times\frac{(2n+\abs{\Gam}+d+2)}{2(2k+\abs{\tilde{\Gam}}+d+1)}\norm{\PP_{n-k}^{\left(2k+\abs{\tilde{\Gam}}+d+1,\gamma_{d+2}\right)}}_{\J\left(2k+\abs{\tilde{\Gam}}+d+1,\gamma_{d+2}\right)}^{2}.
\end{multline}
From the induction hypothesis we have, for each $k \in [n]_0$,
\begin{equation}\label{eq:inductive-hypp}
\norm{q_{k}}_{\tilde{\Gam}}^{2} \leq \frac{(\gamma_{1}+\gamma_{2}+3)(2k+\abs{\tilde{\Gam}}+d+1)}{(\gamma_{1}+1)(\gamma_{2}+1)}\norm{q_{k}}_{\mathrm{ridge},\tilde{\Gam}+e_{1}+e_{2}}^{2}.
\end{equation}
Combining \eqref{eq:normrn}, \eqref{eq:identnormT} and \eqref{eq:inductive-hypp}, using \autoref{prop:normT} (with shifted $d$, $\gamma_1$ and $\gamma_2$) to make squared $\W_{\mathrm{ridge},\Gam+e_1+e_2}$-norms of the $\mathcal{T}_{n,k}^{\Gam+e_{1}+e_{2}}(q_{k})$ appear and using \autoref{descompobreveV} to combine the latter into the corresponding squared norm of $p$, we obtain
\begin{equation*}
\norm{p}_{\Gam}^{2} \leq \frac{(\gamma_{1}+\gamma_{2}+3)(2n+\abs{\Gam}+d+2)}{(\gamma_{1}+1)(\gamma_{2}+1)}\norm{p}_{\mathrm{ridge},\Gam+e_{1}+e_{2}}^{2},
\end{equation*}
which is the desired result in dimension $d+1$ and completes the induction argument.
\end{proof}
\end{lemma}
\begin{thm}\label{thm:whole-depleted-Schur}
Let $d\in\mathbb{N}_{\geq 2}$, $\Gam\in(-1,\infty)^{d+1}$ and $n\in\mathbb{N}_{0}$.
Then, for all $i, j \in [d+1]$ with $i \neq j$ and all polynomials $p \in \Pi_{n}^d$,
\begin{equation*}
\norm{p}_{\Gam}^{2} \leq \frac{(\gamma_{i}+\gamma_{j}+3)(n+1)(n+\abs{\Gam}+d+1)}{(\gamma_{i}+1)(\gamma_{j}+1)}\norm{\frac{p}{\sqrt{x_{i}+x_{j}}}}_{\Gam+e_{i}+e_{j}}^{2}.
\end{equation*}
\begin{proof}
Let $p\in\Pi_{n}^{d}$.
From the decomposition $\Pi_{n}^{d} = \bigoplus_{k=0}^{n}\mathcal{V}_{k}^{\mathrm{ridge},\Gam}$ stated in \eqref{eq:dep-orthogonal-decomposition-of-polynomials}, the triangle inequality and \autoref{lem:slice-depleted-Schur} we have
\begin{equation*}
\begin{split}
\norm{p}_{\Gam} & \leq \sum_{j=0}^{n}\norm{\proj_{j}^{\mathrm{ridge},\Gam+e_{1}+e_{2}}(p)}_{\Gam}\\
& \leq \sum_{j=0}^n \sqrt{\frac{(\gamma_{1}+\gamma_{2}+3)(2j+\abs{\Gam}+d+1)}{(\gamma_{1}+1)(\gamma_{2}+1)}} \norm{\proj_{j}^{\mathrm{ridge},\Gam+e_{1}+e_{2}}(p)}_{\mathrm{ridge},\Gam+e_{1}+e_{2}}.
\end{split}
\end{equation*}
Hence, via the Cauchy--Schwarz inequality in $\mathbb{R}^{n+1}$ and Parseval's identity, we obtain
\begin{equation*}
\begin{split}
\norm{p}_{\Gam}^{2} & \leq \frac{(\gamma_{1}+\gamma_{2}+3)}{(\gamma_{1}+1)(\gamma_{2}+1)}\sum_{j=0}^{n}(2j+\abs{\Gam}+d+1)\norm{p}_{\mathrm{ridge},\Gam+e_{1}+e_{2}}^{2}\\
& = \frac{(\gamma_{1}+\gamma_{2}+3)(n+1)(n+\abs{\Gam}+d+1)}{(\gamma_{1}+1)(\gamma_{2}+1)} \norm{p}_{\mathrm{ridge},\Gam+e_{1}+e_{2}}^{2}.
\end{split}
\end{equation*}
Recalling the definition of the weight $\W_{\mathrm{ridge},\Gam+e_1+e_2}$ in \eqref{eq:depleted-weight}, we recognize the desired estimate of the case $i = 1$ and $j = 2$.

Next, given $i,j\in[d+1]$ with $i\neq j$, let $\tau\in\mathcal{S}_{d+1}$ be any permutation satisfying $\tau(1) = i$ and $\tau(2) = j$.
Since $T_{\tau}^{\star}p\in\Pi_{n}^{d}$, we can apply the previous estimate to obtain
\begin{equation*}
\norm{T_{\tau}^{\star}p}_{\tau(\Gam)}^{2} \leq \frac{(\gamma_{i}+\gamma_{j}+3)}{(\gamma_{i}+1)(\gamma_{j}+1)}(n+1)(n+\abs{\Gam}+d+1)\norm{T_{\tau}^{\star}p}_{\mathrm{ridge},\tau(\Gam)+e_{1}+e_{2}}^{2}.
\end{equation*}
Now, from the definitions of our affine maps in \eqref{eq:affine-map-T-sigma} and of their associated pullbacks in \eqref{eq:pullback-T-sigma} (recall also that in the introduction we set $x_{d+1} = 1 - \abs{x}$) we have that
\begin{equation*}
\frac{T_{\tau}^{\star}p(x)}{\sqrt{x_{1}+x_{2}}} = \frac{T_{\tau}^{\star}p(x)}{T_{\tau}^{\star}\circ T_{\tau^{-1}}^{\star}(\sqrt{x_{1}+x_{2}})} = \frac{T_{\tau}^{\star}p(x)}{T_{\tau}^{\star}\sqrt{x_{i}+x_{j}}} = T_{\tau}^{\star}\left(\frac{p(x)}{\sqrt{x_{i}+x_{j}}}\right).
\end{equation*}
As $p$ is a polynomial and $\frac{x_i x_j}{x_i+x_j} \leq 1$, the function $\frac{p(x)}{\sqrt{x_i+x_j}}$ belongs to $\LL^2_{\Gam+e_i+e_j}$.
Hence, by part~\ref{prop6onehalf} of \autoref{prop6} and the above chain of equalities,
\begin{equation*}
\norm{\frac{p}{\sqrt{x_{i}+x_{j}}}}_{\Gam+e_{i}+e_{j}}^{2}
= \norm{T_{\tau}^{\star}\left(\frac{p}{\sqrt{x_{i}+x_{j}}}\right)}_{\tau(\Gam+e_{i}+e_{j})}^{2}
= \norm{T_{\tau}^{\star}p}_{\mathrm{ridge},\tau(\Gam)+e_{1}+e_{2}}^{2}.
\end{equation*}
Also, again from part~\ref{prop6onehalf} of \autoref{prop6},
\begin{equation*}
\norm{T_{\tau}^{\star}p}_{\tau(\Gam)}^{2} = \norm{p}_{\Gam}^{2}.
\end{equation*}
The desired result follows.
\end{proof}
\end{thm}

\section{Markov Inequality}\label{sec:markovinequality}
In this section we prove the Markov-type inequality \eqref{crude-Markov}.
With \autoref{thm:whole-depleted-Schur} at hand, our next goal is to combine this estimate with a Bernstein-type inequality derived from the results in \cite{GeXu:2024}.

\begin{proposition}\label{prebernstein}
Let $d\in\mathbb{N}_{\geq2}$, $\Gam\in(-1,\infty)^{d+1}$ and $n\in\mathbb{N}_{0}$.
Then, for all $p\in\Pi_{n}^{d}$ and $x\in\tri^{d}$ there holds:
\begin{equation*}
\left[ \partial_{1,2}p(x) \right]^{2}\frac{x_{1}x_{2}}{x_{1}+x_{2}} \leq (x\cdot \nabla p(x))^{2}\frac{1-\abs{x}}{\abs{x}} + \sum_{1\leq i< j\leq d}\left[ \partial_{i,j}p(x)\right]^{2}\frac{x_{i}x_{j}}{\abs{x}}
\end{equation*}
\begin{proof}
Given $x\in\tri^{d}$, we define the matrices $K(x),M(x)\in\mathbb{R}^{d\times d}$ by
\begin{equation*}
K(x) := \frac{x_{1}x_{2}}{x_{1}+x_{2}}(e_{1}-e_{2})\otimes (e_{1}-e_{2}) = \frac{x_{1}x_{2}}{x_{1}+x_{2}}\begin{bmatrix}
1 & -1 & 0 & \cdots & 0\\
-1 & 1 & 0 & \cdots & 0\\
0 & 0 & 0 & \cdots & 0\\
\vdots & \vdots &\vdots & \ddots & \vdots\\
0 & 0 & 0 & \cdots & 0\\
\end{bmatrix}
\end{equation*}
and
\begin{equation*}
M(x) := \frac{1-\abs{x}}{\abs{x}}x\otimes x + \sum_{1\leq i<j\leq d} \frac{x_{i}x_{j}}{\abs{x}}(e_{i}-e_{j})\otimes (e_{i}-e_{j}).
\end{equation*}
It is easy to check that $K(x)$ and $M(x)$ are symmetric matrices.
We also claim that $M(x)$ is positive-definite; indeed, as $x_{1},\dots,x_{d}>0$ and $1-\abs{x}>0$, for all $z\in\mathbb{R}^{d}$,
\begin{equation*}
z\cdot M(x)z = \frac{1-\abs{x}}{\abs{x}}(x\cdot z)^{2} + \sum_{1\leq i<j\leq d} \frac{x_{i}x_{j}}{\abs{x}}(z_{i}-z_{j})^{2}\geq 0
\end{equation*}
and if $z\cdot M(x)z = 0$, then $z = 0$.

Now let us note that the vector $v^{(1)}(x) := -x_2 e_1 + x_1 e_2$ is an eigenvector of the generalized eigenvalue problem: Find $\lambda \in \mathbb{R}$ and $v \in \mathbb{R}^{d} \setminus \{0\}$ such that
\begin{equation}\label{GEV}
K(x) v = \lambda M(x) v.
\end{equation}
Indeed, it is straightforward to check that
\begin{equation*}
K(x) v^{(1)}(x) = M(x) v^{(1)}(x)
\end{equation*}
which goes on to show that the eigenvalue associated with the eigenvector $v^{(1)}(x)$ is $\lambda^{(1)}(x) = 1$.
As $K(x)$ is symmetric and $M(x)$ is symmetric and positive-definite, we can complete the set $\{(\lambda^{(1)}(x), v^{(1)}(x))\}$ into a complete set $\{(\lambda^{(i)}(x), v^{(i)}(x))\}_{i=1}^{d}$ of eigenpairs of \eqref{GEV}, so that $\{ v^{(i)}(x) \}_{i=1}^{d}$ is an $M(x)$-orthogonal basis of $\mathbb{R}^{d}$.
As $K(x)$ is a rank-1 matrix and $\lambda^{(1)}(x) \neq 0$, it follows that $\lambda^{(i)}(x) = 0$ for $i \in \{2,\dots ,d\}$.
Therefore, $\lambda^{(1)}(x) = 1$ is the maximum of the Rayleigh quotient $(v \cdot K(x) v) / (v \cdot M(x) v)$ among the $v \in \mathbb{R}^{d} \setminus \{0\}$, whence
\begin{equation*}
\left(\forall\,v \in \mathbb{R}^{d}\right) \quad v \cdot K(x) v \leq v \cdot M(x) v.
\end{equation*}
As
\begin{equation*}
\nabla p(x) \cdot K(x) \nabla p(x) = [\partial_{1,2} p(x)]^2 \frac{x_1 \, x_2}{x_1+x_2}
\end{equation*}
and
\begin{equation*}
\nabla p(x) \cdot M(x) \nabla p(x) = (x \cdot \nabla p(x))^2 \, \frac{1-\abs{x}}{\abs{x}}
+ \sum_{1 \leq i < j \leq d} [\partial_{i,j} p(x)]^2 \, \frac{x_i \, x_j}{\abs{x}},
\end{equation*}
we obtain the desired result.
\end{proof}
\end{proposition}
We now cite a weak form of the Sturm--Liouville operator $\mathcal{L}^{\Gam}$ of \eqref{SL-O-strong}.
\begin{proposition}\label{propGeXu}
Let $d\in\mathbb{N}_{\geq2}$, $\Gam\in(-1,\infty)^{d+1}$ and $f,g\in \mathrm{C}^{2}(\tri^{d})$.
Then,
\begin{multline*}
\int_{\tri^{d}}\mathcal{L}^{\Gam}f(x)\cdot g(x)\W_{\Gam}(x)\ \textup{d}x = \int_{\tri^{d}}(x\cdot \nabla f(x))(x\cdot \nabla g(x))(1-\abs{x})\W_{\Gam}(x)\frac{\textup{d}x}{\abs{x}}\\ + \sum_{1\leq i<j\leq d}\int_{\tri^{d}}\partial_{i,j}f(x)\partial_{i,j}g(x)x_{i}x_{j}\W_{\Gam}(x)\frac{\textup{d}x}{\abs{x}}.
\end{multline*}
\begin{proof} This is \cite[Th.~2.4]{GeXu:2024}.
\end{proof}
\end{proposition}

The following Bernstein-type inequality is a variant of \cite[Eq.~(3.2)]{GeXu:2024} with explicit constants.

\begin{proposition}\label{prop:Bernstein12}
Let $d\in\mathbb{N}_{\geq 2}$, $\Gam \in(-1,\infty)^{d+1}$ and $n\in\mathbb{N}_{0}$.
Then,
\begin{equation*}
(\forall\,p\in\Pi_{n}^{d})\quad\norm{\partial_{1,2}p}_{\mathrm{ridge},\Gam+e_{1}+e_{2}}\leq \sqrt{n(n+\abs{\Gam}+d)} \norm{p}_{\Gam}.
\end{equation*}
\begin{proof}
Let $p\in\Pi_{n}^{d}$.
By \autoref{prebernstein}, we deduce the inequality
\begin{multline*}
\norm{\partial_{1,2}p}_{\mathrm{ridge},\Gam+e_{1}+e_{2}}^{2} \leq \int_{\tri^{d}}\left(x\cdot\nabla p\right)^{2}(1-\abs{x})\W_{\Gam}(x)\frac{\textup{d}x}{\abs{x}}\\
+ \sum_{1\leq i < j\leq d}\int_{\tri^{d}}\left[ \partial_{i,j}p(x) \right]^{2}x_{i}x_{j}\W_{\Gam}(x)\frac{\textup{d}x}{\abs{x}}.
\end{multline*}
Furthermore, applying \autoref{propGeXu} to the right hand side yields the bound
\begin{equation}\label{eq:pre-bernstein}
\norm{\partial_{1,2}p}_{\mathrm{ridge},\Gam+e_{1}+e_{2}}^{2}\leq \abs{\langle
\mathcal{L}^{\Gam}p,p \rangle_{\Gam}}\leq \norm{\mathcal{L}^{\Gam}p}_{\Gam}\norm{p}_{\Gam}.
\end{equation}
On the other hand,
\begin{multline}\label{eq:L-bound}
\norm{\mathcal{L}^{\Gam}p}_{\Gam}
= \norm{\sum_{k=0}^n \mathcal{L}^{\Gam} \proj^{\Gam}_k(p)}_{\Gam} \stackrel{\eqref{SL-O-strong}}{=} \norm{\sum_{k=0}^n k(k+\abs{\Gam}+d) \proj^{\Gam}_k(p)}_{\Gam}\\
\stackrel{\text{Parseval}}{\leq} n(n+\abs{\Gam}+d) \norm{p}_{\Gam}.
\end{multline}
Then, combining \eqref{eq:pre-bernstein} with \eqref{eq:L-bound} we have the desired inequality.
\end{proof}
\end{proposition}

\begin{thm}\label{thm:Markovineq}
Let $d\in\mathbb{N}_{\geq 2}$, $\Gam\in(-1,\infty)^{d+1}$ and $n\in\mathbb{N}_{0}$.
Then,
\begin{equation*}
\left(\forall\,p\in\Pi_{n}^{d}\right)\quad \norm{\nabla p}_{\Gam} \leq \sqrt{\frac{1}{\gamma_{d+1}+1} \sum_{i=1}^{d}\frac{\gamma_{i}+\gamma_{d+1}+3}{\gamma_{i}+1}} \, n (n+\abs{\Gam}+d)\norm{p}_{\Gam}.
\end{equation*}
\begin{proof}
Let $p \in \Pi_{n}^{d}$.
First, we observe that combining \autoref{thm:whole-depleted-Schur} with \autoref{prop:Bernstein12} we have the estimate
\begin{multline}\label{eq:Markov12}
\norm{\partial_{1,2}p}_{\Gam}^{2}\leq \frac{(\gamma_{1}+\gamma_{2}+3)}{(\gamma_{1}+1)(\gamma_{2}+1)}n(n+\abs{\Gam}+d)\norm{\partial_{1,2}p}_{\mathrm{ridge},\Gam+e_{1}+e_{2}}^{2}\\
\leq \frac{(\gamma_{1}+\gamma_{2}+3)}{(\gamma_{1}+1)(\gamma_{2}+1)}n^{2}(n+\abs{\Gam}+d)^{2}\norm{p}_{\Gam}^{2}.
\end{multline}
Now, since $T_{\sigma_{1}}^{\star}p \in \Pi_{n}^{d}$, it follows from \eqref{eq:Markov12} and the definition of the permutation $\sigma_1$ at \autoref{def:Sigmad} that
\begin{equation*}
\norm{\partial_{1,2}T_{\sigma_{1}}^{\star}p}_{\sigma_{1}(\Gam)}^{2} \leq \frac{(\gamma_{2}+\gamma_{d+1}+3)}{(\gamma_{2}+1)(\gamma_{d+1}+1)}n^{2}(n+\abs{\Gam}+d)^{2} \norm{T_{\sigma_{1}}^{\star}p}_{\sigma_{1}(\Gam)}^{2}.
\end{equation*}
Applying part~\ref{prop8one} in \autoref{prop8} and part~\ref{prop6onehalf} of \autoref{prop6} to the above inequality results in
\begin{equation*}
\norm{\partial_{2}p}_{\Gam}^{2} \leq \frac{(\gamma_{2}+\gamma_{d+1}+3)}{(\gamma_{2}+1)(\gamma_{d+1}+1)}n^{2}(n+\abs{\Gam}+d)^{2} \norm{p}_{\Gam}^{2}.
\end{equation*}
By analogous arguments, we obtain the corresponding bounds for $\norm{\partial_{i}p}_{\Gam}^{2}$ for the rest of the $i \in [d]$, and so,
\begin{equation*}
(\forall\,i\in[d]) \quad \norm{\partial_{i}p}_{\Gam}^{2} \leq \frac{(\gamma_{i}+\gamma_{d+1}+3)}{(\gamma_{i}+1)(\gamma_{d+1}+1)}n^{2}(n+\abs{\Gam}+d)^{2} \norm{p}_{\Gam}^{2}.
\end{equation*}
Then, the desired estimate follows by summing the above bound over $i \in [d]$.
\end{proof}
\end{thm}

\begin{proposition}\label{pro:MarkovIsSharp}
Let $d\in\mathbb{N}_{\geq2}$, $n\in\mathbb{N}_0$ and $\Gam\in(-1,\infty)^{d+1}$.
Then,
\begin{multline*}
2\sqrt{d\mathfrak{R}^{(\abs{\tilde{\Gam}}+d-1,\gamma_{d+1})}(n)} \leq \sup_{p\in\Pi_{n}^{d}-\{0\}}\frac{\norm{\nabla p}_{\Gam}}{\norm{p}_{\Gam}} \\
\leq \sqrt{\frac{1}{\gamma_{d+1}+1}\sum_{i=1}^{d}\frac{\gamma_{i}+\gamma_{d+1}+3}{\gamma_{i}+1}}n (n+\abs{\Gam}+d)
\end{multline*}
Moreover, the power of $n$ in \autoref{thm:Markovineq} is sharp, that is, if there exists $C>0$, $r\in\mathbb{R}$ and $n_{0}\in\mathbb{N}_{0}$ such that
\begin{equation*}
(\forall n\geq n_{0})\ (\forall p \in\Pi_{n}^{d})\quad \norm{\nabla p}_{\Gam} \leq Cn^{r}\norm{p}_{\Gam}.
\end{equation*}
Then, $r\geq 2$.
\begin{proof}
The upper bound is given by \autoref{thm:Markovineq}.
For the lower bound, let $Q\in\Pi_{n}^{d}$ such that $Q(x) = \PP_{n}^{(\abs{\tilde{\Gam}} + d,\gamma_{d+1})}(1-2\abs{x})$.
It is straightforward to check that $\partial_{i}Q_{n} = -2\left(\PP_{n}^{(\abs{\tilde{\Gam}} + d,\gamma_{d+1})}\right)'(1-2\abs{x})$ for all $i\in[d]$.
Then, performing the change of variables $x = \Psi_{d}(\zeta,\eta)$ (cf.~\eqref{Psi} and \eqref{Psi-d}) on the integrals underlying the $\W_{\Gam}$-weighted norms of $Q_{n}$ and $\partial_{i}Q_{n}$, we obtain
\begin{equation*}
\frac{\norm{\partial_{i}Q_{n}}_{\Gam}^{2}}{\norm{Q_{n}}_{\Gam}^{2}} = 4\frac{\norm{\left(\PP_{n}^{(\abs{\tilde{\Gam}} + d,\gamma_{d+1})}\right)'}_{\J(\abs{\tilde{\Gam}}+d-1,\gamma_{d+1})}^{2}}{\norm{\PP_{n}^{(\abs{\tilde{\Gam}} + d,\gamma_{d+1})}}_{\J(\abs{\tilde{\Gam}}+d-1,\gamma_{d+1})}^{2}} \overset{\text{\autoref{prop:Markov-sharpness-1d}}}{=} 4\mathfrak{R}^{(\abs{\tilde{\Gam}}+d-1,\gamma_{d+1})}(n).
\end{equation*}
Thus,
\begin{equation*}
2\sqrt{d\mathfrak{R}^{(\abs{\tilde{\Gam}}+d-1,\gamma_{d+1})}(n)} = \frac{\norm{\nabla Q_{n}}_{\Gam}}{\norm{Q_{n}}_{\Gam}} \leq \sup_{p\in\Pi_{n}^{d}}\frac{\norm{\nabla p}_{\Gam}}{\norm{p}_{\Gam}}
\end{equation*}
The above inequality imply that, there exists a positive constant $C>0$ such that
\begin{equation*}
Cn^{2}\leq\sup_{p\in\Pi_{n}^{d}}\frac{\norm{\nabla p}_{\Gam}}{n^{2}\norm{p}_{\Gam}}
\end{equation*}
Hence, the sharpness of the constant in the Markov-type inequality holds.
\end{proof}
\end{proposition}

For future reference, we encapsulate an iterated version of \autoref{thm:Markovineq}.

\begin{corollary}\label{coromarkov}
Given $d\in\mathbb{N}_{\geq2}$, let $\Gam\in(-1,\infty)^{d+1}$, $r\in\mathbb{N}$ and $n\in\mathbb{N}_{0}$.
Then, there exists a constant $C >0$, depending only on $d$, $r$ and $\Gam$, such that
\begin{equation*}
(\forall p\in\Pi_{n}^{d})~\norm{\nabla^{r}p}_{\Gam}\leq C n^{2r}\norm{p}_{\Gam}.
\end{equation*}
\begin{proof}
We proceed by induction on the order of differentiation, observing that the base case corresponds to \autoref{thm:Markovineq}.
\end{proof}
\end{corollary}

\bibliographystyle{amsalpha}
\bibliography{triangle-SOP-refs}

\end{document}